\documentclass[11pt]{article}

\usepackage[ansinew, utf8]{inputenc}
\usepackage{mathtools, amsthm, amssymb, bbm}
\usepackage{multicol,latexsym, array, xcolor}
\usepackage{paralist}
\usepackage{dsfont}
\usepackage{color}
\usepackage{mathrsfs}
\usepackage{natbib}

\usepackage{bm}
\usepackage[colorlinks,citecolor=blue,urlcolor=blue]{hyperref}
\usepackage{graphicx} 
\usepackage[rightcaption]{sidecap}
\usepackage{cleveref}
\usepackage[font=footnotesize,labelfont=bf]{caption}
\usepackage{subcaption}
\usepackage{wrapfig}
\usepackage{nicefrac}
\usepackage{tikz}

\numberwithin{equation}{section}

\usepackage{wrapfig}

	\DeclareMathOperator*{\N}{\mathbb{N}}  
	\DeclareMathOperator*{\R}{\mathbb{R}}
	
	\DeclareMathOperator*{\argmin}{arg\,min}

\providecommand{\keywords}[1]{\textbf{\textit{Keywords}} #1}
\providecommand{\keywordsMSC}[1]{\textbf{\textit{MSC 2010 subject classification}} #1}

\newtheorem{theorem}{Theorem}[section]
\newtheorem{definition}[theorem]{Definition}
\newtheorem{example}[theorem]{Example}
\newtheorem{corollary}[theorem]{Corollary}

\newtheorem{remark}[theorem]{Remark}

\newcommand{\footremember}[2]{%
	\footnote{#2}
	\newcounter{#1}
	\setcounter{#1}{\value{footnote}}%
}
\newcommand{\footrecall}[1]{%
	\footnotemark[\value{#1}]%
}

\makeatletter				
\def\namedlabel#1#2{\begingroup
    #2%
    \def\@currentlabel{#2}%
    \phantomsection\label{#1}\endgroup
}
\makeatother

\begin{document}

\author{Marcel Klatt \footremember{ims}{\scriptsize Institute for Mathematical
		Stochastics, University of G\"ottingen,
		Goldschmidtstra{\ss}e 7, 37077 G\"ottingen} 
	\and 
	Carla Tameling \footrecall{ims}{}
	\and 
	Axel Munk \footrecall{ims} \footnote{\scriptsize Max Planck Institute for Biophysical
		Chemistry, Am Fa{\ss}berg 11, 37077 G\"ottingen}}

	\title{Empirical Regularized Optimal Transport: Statistical Theory and Applications}

\maketitle
\begin{abstract}
We derive limit distributions for empirical regularized optimal transport distances between probability distributions supported on a finite metric space and show consistency of the (naive) bootstrap. In particular, we prove that the empirical regularized transport plan itself asymptotically follows a Gaussian law. The theory includes the Boltzmann-Shannon entropy regularization and hence a limit law for the widely applied Sinkhorn divergence. \\
Our approach is based on an application of the implicit function theorem to necessary and sufficient optimality conditions for the regularized transport problem. The asymptotic results are investigated in Monte Carlo simulations. We further discuss computational and statistical applications, e.g. confidence bands for colocalization analysis of protein interaction networks based on regularized optimal transport.
\end{abstract}


\keywordsMSC{Primary: 62E20, 62G20, 65C60 Secondary: 90C25, 90C31, 90C59}

\section{Introduction}\label{Introduction}

The theory of optimal transport (OT) has a long history in physics, mathematics, economics and related areas, see e.g. \citet{monge1781memoire}, \citet{kantorovich1942translocation}, \citet{rachev1998mass}, \citet{villani2008optimal} and \citet{galichon2016optimal}. Recently, also OT based \textit{data analysis} has become popular in many areas of application, among others, in computer science \citep{schmitz2018wasserstein,balikas2018cross}, mathematical imaging \citep{rubner2000earth,ferradans2014regularized,adler2017learning} and machine learning \citep{arjovsky2017wasserstein,lu2017optimal,sommerfeld2018optimal}. \\
In this paper, we are concerned with certain statistical aspects of \textit{regularized} optimal transport (ROT) on a finite number of support points, e.g. representing spatial locations. To this end, let the ground space $\mathcal{X}=\{x_1,\ldots,x_N\}$ be finite and equipped with a metric $d\colon \mathcal{X}\times \mathcal{X}\to \R_{\geq 0}$, where $\R_{\geq 0}$ denotes the non-negative and $\R_{>0}$ ($\R_{<0}$) the positive (negative) reals. Each probability distribution on $\mathcal{X}$ is represented as an element in $\Delta_N$ the $N$-dimensional simplex of vectors $r\in \mathbb{R}^N$ such that $\sum_{i=1}^N r_i=1$, $r_i\geq 0$. For the sake of exposition, we implicitly assume that $r,s\in \Delta_N$ are vectors of the same length. However, this can be easily generalized to different lengths (\Cref{rem:generalsensitivity}). We represent the cost to transport one unit from $x_i$ to $x_j$ as a vector $c_p\in \R^{N^2}$ defined by the underlying metric $d$ with entries ${c_p}_{(i-1)N+j}\coloneqq d^p(x_i,x_j)$ for $p\geq 1$. Determining the OT between the probability distributions $r,s \in \Delta_N$ on $\mathcal{X}$ then amounts to solve the standard linear program
\begin{equation}\label{eq:OT}
\begin{aligned}
\min_{\pi \in \mathbb{R}^{N^2}}&\quad   \langle c_p, \pi \rangle \\
\textup{subject to }&\quad A \pi =  \begin{bmatrix}r\\s\end{bmatrix}\, ,
                   \pi \geq  0\, . \\ 
\end{aligned}
\end{equation}
The coefficient matrix 
\begin{equation*}
A\coloneqq \begin{bmatrix}I_{N\times N} \otimes {\mathbbm{1}_{1\times N}} \\{\mathbbm{1}_{1\times N}} \otimes I_{N \times N} \end{bmatrix} \in \mathbb{R}^{2N\times N^2}\, ,
\end{equation*}
where $\otimes$ denotes the Kronecker product, encodes the marginal constraints for any transport plan $\pi \in \R^{N^2}$ so that considered as a $N\times N$ matrix its row and column sums are equal to $r$ and $s$, respectively. In total, the linear program \eqref{eq:OT} consists of $2N$ linear equality constraints in $N^2$ unknowns.
An optimal solution of \eqref{eq:OT} denoted as $\pi_p(r,s)$ (not necessarily unique) is known as an \textit{optimal transport plan}. The quantity 
\begin{equation}\label{eq:Wasserstein}
W_p(r,s)\coloneqq \left\langle c_p, \pi_p(r,s) \right\rangle^{\frac{1}{p}}\, ,
\end{equation}
that is the $p$-th root of the optimal value of \eqref{eq:OT}, is referred to as the OT distance, Kanto\-rovich-Rubinstein distance, earth mover's distance or $p$-th Wasserstein distance between the probability distributions $r$ and $s$. \\
Despite its conceptual appeal and its first practical success in various (statistical) applications \citep{munk1998nonparametric,del1999tests,evans2012phylogenetic,tameling2017empirical,chernozhukov2017monge,sommerfeld2018inference,panaretos2018statistical}, the routine use of OT based data analysis is still lacking as for many real world applications it is severely hindered by its computational burden to solve \eqref{eq:OT}. The zoo of OT solvers is diverse and classical examples include the Hungarian method \citep{kuhn1955hungarian}, the auction algorithm \citep{bertsekas1989auction} or the transportation simplex \citep{luenberger2008linear}. However, their respective (average) runtime yields a severe limitation for real world instances. For example, the auction algorithm is known to require $\mathcal{O}(N^3\log(N))$ elementary operations.\\
There has been made certain progress to overcome this numerical obstacle. For instance, exploiting properties of the $l_1$-distance as a specific instance for $d$, on a regular grid the OT problem \eqref{eq:OT} can be stated in only $\mathcal{O}(N)$ unknowns as it suffices to consider transport between neighbouring grid points. This results in an algorithm with average time complexity $\mathcal{O}(N^2)$ \citep{ling2007efficient}. For more general distances \citet{gottschlich2014shortlist} introduce the shortlist method, which has been shown empirically to have a runtime of magnitude $\mathcal{O}(N^{2.5})$. Furthermore, multiscale approaches \citep{gerber2017multiscale} and sparse approximate methods have recently been developed, including \citet{schmitzer2016sparse} who computes OT via a sequence of sparse problems solved by an arbitrary exact transportation algorithm. Empirical simulations demonstrate that algorithms, e.g. the transportation simplex, do benefit when applied in a multiscale fashion in terms of memory demand and runtime. However, solving already moderately sized problems in reasonable time is still a challenging issue, e.g. in two or three-dimensional imaging, where $N$ is of magnitude $\sim 10^6$-$10^8$ \citep{schrieber2017dotmark}, and large scale problems, such as temporal-spatial image or network analysis seem currently out of reach.\\
This encouraged the development of various surrogates for the OT distance which are computationally better accessible. We mention thresholding of the full distance leading to graph sparsification \citep{pele2009fast}, relaxation \citep{ferradans2014regularized} and regularized OT distances \citep{dessein2018regularized,essid2018quadratically}. Among the most prominent proposals for the latter approach is the entropy regularization of OT \citep{cuturi2013sinkhorn,peyre2017computational}. Instead of solving the linear program \eqref{eq:OT}, the entropy regularization approach asks to solve for a positive regularization parameter $\lambda>0$ the \textit{regularized} OT (ROT) problem
\begin{equation}\label{eq:EntropyROT}
\begin{aligned}
\min_{\pi \in \mathbb{R}^{N^2}}&\quad  \langle c_p, \pi \rangle + \lambda f(\pi) \\
\textup{subject to}&\quad A \pi =  \begin{bmatrix}r\\s\end{bmatrix}\, .\\
\end{aligned}
\end{equation}
The function $f\colon \R^{N^2} \to \R$ is the negative Boltzmann-Shannon entropy defined for $\pi \in \R^{N^2}$ as 
\begin{equation}\label{eq:entropy}
f(\pi) \coloneqq
\begin{cases}
 \sum_{i=1}^{N^2} \pi_{i}\log(\pi_i)-\pi_i+1 &\text{for }\pi \in \mathbb{R}_{\geq 0}^{N^2}\, ,  \\
+ \infty &\text{otherwise}
\end{cases}
\end{equation}
with the convention $0\log(0)=0$. Different to \eqref{eq:OT}, the regularization in \eqref{eq:EntropyROT} avoids the non-negativity constraint $\pi \geq 0$ since introducing entropy in the objective already enforces any feasible solution to be non-negative. Moreover, ROT \eqref{eq:EntropyROT} is a strictly convex optimization program and hence has a \textit{unique} optimal solution $\pi_{p,\lambda, f}(r,s)$ denoted as \textit{entropy regularized optimal transport plan}. The entropy ROT distance, also known as $p$-th Sinkhorn divergence \citep{cuturi2013sinkhorn}, is then defined as
	\begin{equation}
	\label{eq:sinkhorn}
	W_{p,\lambda,f}(r,s) \coloneqq \left\langle c_p, \pi_{p,\lambda, f}(r,s) \right\rangle^{\frac{1}{p}}\, .
	\end{equation}
The major benefit of entropy regularization is of algorithmic nature. The entropy ROT plan can be approximated by the Sinkhorn-Knopp algorithm originally introduced by \citet{sinkhorn1964relationship} that has a linear convergence rate and only requires $\mathcal{O}(N^2)$ operations in each step \citep{cuturi2013sinkhorn,altschuler2017near}. It has been argued that as the regularization parameter $\lambda>0$ decreases to zero in \eqref{eq:EntropyROT} this approximates the solution of the OT problem and hence serves as a good proxy. Nevertheless, high accuracy is computationally hindered in small regularization regimes \citep{benamou2015iterative,schmitzer2016stabilized} as the runtime of these algorithms scale with $\lambda^{-2}$ \citep{dvurechensky2018computational}.\\
The results of this paper complement these \textit{computational} findings for ROT, as we will show a substantially different \textit{statistical} behaviour of the regularized ($\lambda>0$) compared to non-regularized OT $(\lambda=0$) when the probability distributions $r$ and $s$ are estimated from data, hence randomly perturbed.
To this end and as often typical in applications, the underlying population distribution $r$ (resp. $s$ or both) is estimated from given data by its empirical version 
\begin{equation}\label{eq:empiricalmeasures}
\hat{r}_n=\frac{1}{n}\sum_{i=1}^n \delta_{X_i}, \, \left(\hat{s}_m=\frac{1}{m}\sum_{i=1}^m \delta_{Y_i}\right) 
\end{equation}
derived by a sample of $\mathcal{X}$-valued random variables $X_1,\ldots,X_n\overset{i.i.d.}{\sim} r$ (resp. $Y_1,\ldots,Y_m\overset{i.i.d.}{\sim} s$). In this paper, we explicitly compute the limit distributions (after proper standardization) of the population version of $\pi_{p,\lambda,f}(r,s)$ and $W_{p,\lambda,f}(r,s)$ based on $\hat{r}_n$ and $\hat{s}_m$. Notably, we make use of the fact that entropy regularization \eqref{eq:EntropyROT} enforces a \textit{dense} structure for its entropy ROT plan, i.e. all its entries are positive. This property, long known and favoured, for instance, in traffic prediction \citep{wilson1969use} facilitates our sensitivity analysis for ROT (\Cref{thm:sensitivity}). Our main result (\Cref{thm:optimalsolutionlimit}) concerns the empirical ROT plan itself. More precisely, for $\lambda>0$ in \eqref{eq:EntropyROT} the random quantity $\pi_{p,\lambda,f}(\hat{r}_n,s)$ converges in distribution ($\overset{D}{\longrightarrow}$) to an $N^2$-dimensional Gaussian law
\begin{equation}\label{eq:optimalsolutionlimit}
\sqrt{n}\left\lbrace \pi_{p,\lambda,f}(\hat{r}_n,s)-\pi_{p,\lambda,f}(r,s) \right\rbrace  \overset{D}{\longrightarrow} \mathcal{N}_{N^2}\left(0,\Sigma_{p,\lambda,f}(r \vert s)\right)\, ,
\end{equation}
as $n \to \infty$. The covariance $\Sigma_{p,\lambda,f}(r \vert s)$ depends on $\lambda$, the Hessian of the regularization function $f$ and the probability distributions $r$ and $s$ (see \eqref{eq:covarianceregtransport}). Furthermore, the empirical ROT distance asymptotically follows a centred Gaussian limit law (\Cref{thm:optimaltransportlimit}), that is
\begin{equation}\label{eq:optimaldistancelimit}
\sqrt{n}\left\lbrace W_{p,\lambda,f} (\hat{r}_n,s)-W_{p,\lambda,f}(r,s) \right\rbrace  \overset{D}{\longrightarrow} \mathcal{N}_1\left(0,\sigma^2_{p,\lambda,f} (r \vert s)\right) \, .
\end{equation}
Our results hold true for $r=s$ and $r\neq s$ as well as for other cost vectors $c$ (\Cref{rem:generalsensitivity}) and are valid for a broad class of regularizers $f$ in \eqref{eq:EntropyROT}. These will be denoted as \textit{proper regularizers} (\Cref{def:properregularizer}) and subsume various regularization methods of OT \citep{dessein2018regularized}. It includes the widely applied entropy ROT in \eqref{eq:entropy} and hence limit distributions for the empirical entropy ROT plan and, as a consequence, for the empirical Sinkhorn divergence. For the latter see also \citet{bigot2017central} who obtained limit distributions in a similar fashion to \eqref{eq:optimaldistancelimit} for the optimal value of \eqref{eq:EntropyROT} with a straightforward application of the technique in \citet{sommerfeld2018inference}. Note that the technique we introduce here allows to treat the ROT plan itself also in notable distinction to $\lambda=0$, where such a result as in \eqref{eq:optimalsolutionlimit} is not known. Our limit theorems for the ROT plan will then be used in \Cref{Bootstrap} and \Cref{Colocalization} to derive several statistical consequences.\\
Our findings highlight a substantial difference to related limit laws for \textit{non-regu\-larized} transport. In fact, our approach is different to the technique used in \citet{sommerfeld2018inference} and fails for (non-regularized) OT \eqref{eq:OT} as the optimal (non-regularized) transport plan is known to be sparse. Further, the optimization problem in \eqref{eq:EntropyROT} is non-linear whereas \eqref{eq:OT} is a linear program. The limit distributions for the empirical ROT distance turn out to be substantially different to those of the empirical (non-regularized) OT distance (see \citet{sommerfeld2018inference}). More specifically, for $r=s$ the empirical (non-regularized) transport distance $W_p(\hat{r}_n,r)$ does not follow asymptotically a Gaussian law, in contrast to \eqref{eq:optimaldistancelimit}, as for $p\geq 1$ it holds that
\begin{equation}\label{eq:finitetransportlimit}
n^{1/2p}\,  W_p(\hat{r}_n,r) \overset{D}{\longrightarrow} \left\lbrace \max_{u \in \Phi^*} \langle G,u \rangle \right\rbrace^{1/p}\, ,
\end{equation}
as $n \to \infty$ \citep{sommerfeld2018inference}. Here, $\Phi^*$ is the set of dual solutions for \eqref{eq:OT} and $G$ is a centred $N$-dimensional Gaussian random vector with covariance matrix 
\begin{equation}\label{eq:sigmaone}
\Sigma(r) \coloneqq
\begin{pmatrix}
r_1(1-r_1) & -r_1 r_2 & \ldots & -r_1 r_N\\
-r_2r_1 & r_2(1-r_2)  & \ldots & -r_2 r_N\\
\vdots & \vdots & \ddots & \vdots \\
-r_N r_1 & -r_N r_2 & \ldots & r_N(1-r_N)
\end{pmatrix}\, .
\end{equation}
This different limit behaviour provides some insight into the above mentioned computational difficulties to approximate the non-regularized OT distance by the Sinkhorn divergence as $\lambda \searrow 0$.\\ 
The outline of this paper is as follows. As a prerequisite for our main methodology and results we provide in \Cref{Sensitivity} all necessary terminology from convex optimization. We derive sensitivity results for ROT plans which might be of interest by itself as they describe the stability of ROT plans when perturbing the boundary conditions given by $r$ and $s$. Our approach is based on an application of the implicit function theorem to necessary and sufficient optimality conditions for \eqref{eq:EntropyROT} with more general regularizers. Parallel to our work, such a sensitivity result was partially obtained by \citet{luise2018differential}. However, their proof is carried out on the dual formulation of \eqref{eq:EntropyROT} and is limited to entropy regularization.\\
\Cref{Distributions} is dedicated to distributional limit results stated in \Cref{thm:optimalsolutionlimit} and \Cref{thm:optimaltransportlimit}. Moreover, we give rates for the regularizer $\lambda(n)$ tending to zero and depending on the sample size in order to asymptotically recover the limit laws for (non-regularized) OT (\Cref{Scaling}). Since our proof for \eqref{eq:optimalsolutionlimit} and \eqref{eq:optimaldistancelimit} is based on a delta method, as a byproduct we obtain consistency of the (naive) $n$ out of $n$ bootstrap in \Cref{Bootstrap}. This is again in notable contrast to the (non-regularized) empirical OT distance, where the (naive) $n$ out of $n$ bootstrap is known to fail \citep{sommerfeld2018inference}. In \Cref{Simulations} we investigate in a Monte Carlo study the approximation of the empirical Sinkhorn divergence sample distribution by its theoretical limit law. In addition, we analyse the influence of the amount of regularization $\lambda$, compare our results to asymptotic distributions for the non-regularized OT distance and investigate the bootstrap empirically.\\
\Cref{Subsampling} introduces a resampling method that reduces the computational complexity still inherent in the computation of ROT, especially for large instance sizes. Based on this, we propose a protocol to analyse the colocalization for images that are beyond the scope of computational feasibility due to their representation by several thousands of pixels. In \Cref{Colocalization} we utilize the ROT plan and our resampling method for a statistical analysis of colocalization for protein interaction networks in cells.\\
Finally, we stress that while the Sinkhorn divergence is numerically appealing, its interpretation is not as straightforward. Although it shares some distance-like properties \citep{cuturi2013sinkhorn}, $W_{p,\lambda,f}$ is not a distance, in particular $W_{p,\lambda,f}(r,r)>0$ when $\lambda>0$. This hinders a simple interpretation, as this value will depend on the specific distribution $r$. Nevertheless, as we show in this paper, regularization allows for a rigorous statistical analysis of the corresponding ROT \textit{plan}, a task which is currently out of sight for the (non-regularized) OT plan. Compared to any OT based distance (\textit{regularized} or not), the ROT \textit{plan} encodes more structural information across scales and hence serves as a more informative tool for inferential statistics. This is utilized for the analysis of protein networks in \Cref{Colocalization} where we provide resampling based confidence bands (\Cref{thm:colocalizationlimit}) for a measure of protein proximity based on the estimated ROT between two protein distributions in a cell compartment. \\
To ease readability the proofs are deferred to the Appendix \ref{Appendix}.

\section{Sensitivity Analysis for Regularized Transport Plans}\label{Sensitivity}

OT in its standard form \eqref{eq:OT} can be stated in terms of only $2N-1$ equality constraints instead of $2N$. In fact, since the total supply equals the total demand $(r,s \in \Delta_N)$ any one of the equality constraints is redundant. Without loss of generality we delete the last constraint. Consequently, we define for $r,s \in \Delta_N$ the feasible set for OT as
\begin{equation}\label{eq:feasibleset}
\mathcal{F}(r,s) \coloneqq \left\lbrace \pi \in \mathbb{R}_{\geq 0}^{N^2}\, \bigr\rvert \, A_{\star}\pi =[r,s_{\star}]^T \right\rbrace
\end{equation}
with coefficient matrix $A_{\star}$ and vector $s_{\star}$, where the subscript star denotes the deletion of the last row of the matrix $A$ in \eqref{eq:OT} and the last entry of the vector $s \in \Delta_N$, respectively. This reduction allows for linearly independent constraints described by $A_{\star}$ or equivalently full rank of $A_{\star}^T$ \citep{luenberger2008linear}.

\begin{remark}\label{rem:asymetric}
In the sequel, we consider OT \eqref{eq:OT} or the ROT \eqref{eq:EntropyROT} as an optimization program with feasible set $\mathcal{F}(r,s)$ involving $A_{\star}$ and $s_{\star} \in (\Delta_N)_\star$. Hence, this requires some caution as we treat $r$ and $s$ asymmetrically. 
\end{remark}

We consider regularizers $f$ in \eqref{eq:EntropyROT} that are of Legendre type, which means that $f$ is a closed proper convex function on $\R^{N^2}$ which is essentially smooth and strictly convex on the interior of its domain \citep{dessein2018regularized}.
Recall that a function is essentially smooth if it is differentiable on the interior of its domain and for every sequence $(x_k)_{k \in \N}\subset \text{int}(\text{dom }f)$ converging to a boundary point of $\text{int}(\text{dom }f)$ it holds that $\lim_{k \to \infty} \Vert \nabla f(x_k) \Vert=+\infty$. For further details on the class of Legendre functions we refer to \citet{rockafellar1970}.

\begin{definition}[Proper Regularizer] \label{def:properregularizer}
Let $f\colon \R^{N^2}\to \R\cup\{+\infty\}$ be twice continuously differentiable on the interior of its domain with positive definite Hessian $\nabla^2 f$. Moreover, assume for $f$ and its Fenchel conjugate $f^*$ that 
\begin{multicols}{2}
\begin{enumerate}
\item $f$ is of Legendre type,
\item $\R_{< 0}^{N^2} \subset \text{dom } f^*$,
\item $(0,1)^{N^2} \subseteq \text{dom }f$,
\item $\text{dom }f \subseteq \mathbb{R}_{\geq 0}^{N^2}$.
\end{enumerate}
\end{multicols}
\noindent Then $f$ is said to be a proper regularizer.
\end{definition}

Examples for proper regularizers are discussed more carefully in the next subsection.
By strict convexity, for each proper regularizer $f$, regularization parameter $\lambda>0$ and $p\geq 1$, we obtain a \textit{unique} ROT plan
\begin{equation}\label{eq:ROTsolution}
\pi_{p,\lambda,f}(r,s)=\argmin_{\pi \in \mathcal{F}(r,s)}\, \langle c_p, \pi \rangle + \lambda f(\pi)\, .
\end{equation}

As the main contribution of this section, we provide a \textit{sensitivity analysis} of \eqref{eq:ROTsolution} in the sense of non-linear programming, i.e. we investigate how the ROT plan $\pi_{p,\lambda,f}(r,s)$ changes when perturbing the marginal constraints given by $r,s \in \Delta_N$. For a general introduction to sensitivity analysis we refer to \citet{fiacco1984sensitivity}.

\begin{theorem}\label{thm:sensitivity}
Let $f$ be a proper regularizer, $r_0,s_0 \in \Delta_N$ and $\lambda>0$.
\begin{enumerate}
\item The regularized transport plan $\pi_{p,\lambda,f}(r_0,s_0)$ in \eqref{eq:ROTsolution} is unique and contained in the positive orthant $\R_{> 0}^{N^2}$.
\item There exists a neighbourhood $\mathcal{U} \subseteq \Delta_N \times (\Delta_N)_\star$ of $(r_0,{s_0}_{\star})$ and a unique continuously differentiable function
\begin{equation*}
\phi_{p,\lambda,f}\colon \mathcal{U} \to \mathbb{R}^{N^2}
\end{equation*}
such that $\phi_{p,\lambda,f}(r_0,{s_0}_{\star})=\pi_{p,\lambda,f}(r_0,s_0)$. Moreover, the regularized transport plan is parametrized by $\phi_{p,\lambda,f}$ for all $(r,s_{\star}) \in \mathcal{U}$ with derivative at $(r_0,{s_0}_{\star})$ given by
\begin{align*}
\nabla \phi_{p,\lambda,f}&(r_0,{s_0}_{\star})= \\
&[\nabla^2 f(\pi_{p,\lambda,f}(r_0,s_0))]^{-1}A_{\star}^T[A_{\star}[\nabla^2 f(\pi_{p,\lambda,f}(r_0,s_0))]^{-1}A_{\star}^T]^{-1}\, ,
\end{align*}
a matrix of size $N^2 \times (2N-1)$.
\end{enumerate}
\end{theorem}

The crucial observation in our proof is that proper regularizers enforce the ROT plan $\pi_{p,\lambda,f}$ to be dense, i.e. contained in the positive orthant. Hence, from an optimization point of view the only binding constraints (constraints fulfilled with equality at the optimal solution) are given by the row and column sum requirement, i.e. $A_{\star}\pi= \begin{bmatrix}r, s_{\star}\end{bmatrix}^T$.
The gradients of these constraints with respect to $\pi$ are linearly independent by full rank of $A_{\star}^T$, a property well-known as the \textit{linear independence constraint qualification}. This allows for an analysis in the spirit of the implicit function theorem. Besides the uniqueness issue in (non-regularized) OT, mimicking the proof for $\lambda=0$ does not work as we require knowledge about the binding constraints additionally inherent in $\pi\geq 0$. The linear independence constraint qualification fails to hold especially in the case where $\pi_p(r,s)$ consists of more than $2N-1$ zeroes, a situation known in linear programming as degeneracy \eqref{eq:OT} (see \citet{luenberger2008linear}).

\begin{remark}[General cost, varying number of support points]\label{rem:generalsensitivity}
\Cref{thm:sensitivity} holds independent of the choice for any non-negative cost vector $c$. However, as we are interested in the ROT distance, we restrict ourselves in the subsequent analysis usually to the case that ${c_p}_{(i-1)N+j}\coloneqq d^p(x_i,x_j)$ for $p\geq 1$. Furthermore, \Cref{thm:sensitivity} implicitly covers different number of support points of the probability distributions $r$ and $s$ which amounts to the situation that, e.g. the probability distribution $r$ does not need to have full support. In such a case, we simply delete the zero entries and our sensitivity result holds for the reduced problem.
\end{remark}

\begin{remark}[Extensions beyond finite spaces]
Our method of proof of \Cref{thm:sensitivity} does not extend to the countable nor to the continuous formulation of ROT. Already for the countable case existence and uniqueness of a ROT plan is not straightforward. For example, the optimal value of ROT \eqref{eq:EntropyROT} can easily be infinity if the marginals are not restricted to have finite $p$-th moment and finite entropy (see \citet{kovavcevic2015entropy}). Note further, that the approach by \citet{tameling2017empirical} for the treatment of OT distance when the ground space is countable is not applicable. This is due to the definition of the Sinkhorn divergence \eqref{eq:sinkhorn} which requires to consider the optimal solution rather than the optimal value of a convex optimization problem.
\end{remark}

\subsection{Proper Regularizers} 

The class of proper regularizers in \Cref{def:properregularizer} is rather rich and some common ones are the negative Boltzmann-Shannon entropy \eqref{eq:entropy} or $l_p$ quasi norms $(0<p<1)$ defined as $f(\pi)= \sum_{i=1}^{N^2} \pi_{i}^p$ for $\pi \in \mathbb{R}_{\geq 0}$ among others. Further examples that have also been the focus of recent research \citep{dessein2018regularized} are given in Table \ref{table:properregularizer}. Moreover, we give their Hessians which are required for the sensitivity analysis (Theorem \ref{thm:sensitivity}).

\renewcommand{\arraystretch}{1.5}
\setlength\tabcolsep{2pt}
\begin{table}[ht]
\begin{center}
\caption{Proper regularizers.}
\fbox{\begin{tabular}{cccc}
Regularizer $f$ & $\text{dom }f$/$f^*$ & $\nabla^2 f$\\ \hline
Boltzmann-Shannon entropy \\ \small $\sum_{i=1}^{N^2} \log(\pi_{i})\pi_{i}-\pi_{i}+1$ & \small $\mathbb{R}_{\geq 0}^{N^2}$/$\mathbb{R}^{N^2}$ & \small diag($\frac{1}{\pi}$) \\[1.4ex]
Burg entropy \\ \small $\sum_{i=1}^{N^2} \pi_{i}-\log(\pi_{i})-1$ & \small $\mathbb{R}_{> 0}^{N^2}$/$(-\infty,1)^{N^2}$ & \small diag($\frac{1}{\pi^2}$) \\[1.4ex]
Fermi-Dirac entropy \\ \small $\sum_{i=1}^{N^2} \log(\pi_{i})\pi_{i}+(1-\pi_{i})\log(1-\pi_{i})$ & \small $[0,1]^{N^2}$/$\mathbb{R}^{N^2}$ & \small diag($\frac{1}{(1-\pi)\pi}$) \\[1.4ex]
$\beta$-potentials ($0<\beta<1$) \\ \small $\frac{1}{\beta(\beta-1)} \sum_{i=1}^{N^2} \pi_{i}^\beta-\beta \pi_{i}+\beta-1$ & \small $\mathbb{R}_{\geq 0}^{N^2}$/$(-\infty,\frac{1}{1-\beta})^{N^2}$ & \small diag($\pi^{\beta-2}$) \\[1.4ex]
$l_p$ quasi norms $(0<p<1)$\\ \small $-\sum_{i=1}^{N^2} \pi_{i}^p$ & \small $\mathbb{R}_{\geq 0}^{N^2}$/$\mathbb{R}_{\leq 0}^{N^2}$ & \small $p(1-p)$\newline diag($\pi^{p-2}$) \\
\end{tabular}}
\label{table:properregularizer}
\end{center}
\end{table}
 
\begin{example}[Sinkhorn divergence]\label{example:entropy}
If $f$ is the negative Boltzmann-Shan\-non entropy in \eqref{eq:entropy}, then the gradient of the parametrization for the entropy ROT plan $\pi_{p,\lambda,f}(r,s)$ is given by 
\begin{equation*}
\nabla \phi_{p,\lambda,f}(r,s_{\star})=D\,A_{\star}^T\,\left[A_{\star}\,D\,A_{\star}^T\right]^{-1}\, ,
\end{equation*}
where $D \in \R^{N^2\times N^2}$ is a diagonal matrix with diagonal $\pi_{p,\lambda,f}(r,s)$.
\end{example}

However, we would like to stress that not all common regularizers fall into the class of proper regularizers.

\begin{example}[A counterexample]
For $f$ the $l_p$ norm $1 \leq p < +\infty$, \Cref{thm:sensitivity} does not hold. In fact, $f$ is not a proper regularizer and allows for a sparse ROT plan \citep{blondel2017smooth}. In particular, we cannot guarantee that its corresponding ROT satisfies the linear independence constraint qualification (see discussion after \Cref{thm:sensitivity}). Hence, our proof strategy for $l_p$ $(p\geq 1)$ ROT plans fails.
\end{example}

\section{Distributional Limits}\label{Distributions}

For two probability distributions $r,s \in \Delta_N$, parameters $\lambda > 0$, $p\geq 1$ and proper regularizer $f$ an estimator for $\pi_{p,\lambda,f}(r,s)$ in \eqref{eq:ROTsolution} is given by its empirical counterpart $\pi_{p,\lambda,f}(\hat{r}_n,s)$ with $\hat{r}_n$ the empirical distribution of the i.i.d. sample $X_1,\ldots,X_n$ in \eqref{eq:empiricalmeasures}. \\
The next theorem states a Gaussian limit distribution for the empirical ROT plan. Since the sensitivity result in \Cref{thm:sensitivity} holds regardless of $r=s$ or $r\neq s$ and as the ROT plan is always dense, we do not derive any substantially difference regarding statistical limit behaviour in either of these cases.

\begin{theorem}\label{thm:optimalsolutionlimit}
Let $r,s \in \Delta_N$ be two probability distributions on the finite metric space $(\mathcal{X},d)$ and let $\hat{r}_n$ be the empirical version given in \eqref{eq:empiricalmeasures} derived by $X_1,\ldots,X_n\overset{i.i.d.}{\sim} r$. 
Then, as the sample size $n$ grows to infinity, it holds that
\begin{equation*}
\sqrt{n}\left\lbrace \pi_{p,\lambda,f}(\hat{r}_n,s)-\pi_{p,\lambda,f}(r,s) \right\rbrace  \overset{D}{\longrightarrow} \mathcal{N}_{N^2}\left(0,\Sigma_{p,\lambda,f}(r \vert s)\right) 
\end{equation*}
with covariance matrix 
\begin{equation}\label{eq:covarianceregtransport}
\Sigma_{p,\lambda,f}(r \vert s) = \nabla_r\, \phi_{p,\lambda,f}(r,s_{\star}) \, \Sigma(r)\, \left[\nabla_r\, \phi_{p,\lambda,f}(r,s_{\star})\right]^T\, ,
\end{equation}
where $\Sigma(r)$ is defined in \eqref{eq:sigmaone} and $\nabla_r\, \phi_{p,\lambda,f}(r,s_{\star})$ are the partial derivatives of $\phi_{p,\lambda,f}$ with respect to $r$ as given in \Cref{thm:sensitivity}.
\end{theorem}

The proof is based on the multivariate delta method and straightforward given \Cref{thm:sensitivity}, hence postponed to the supplement. Further, we prove limit distributions for the empirical counterpart of the ROT distance \eqref{eq:sinkhorn}. Here, $s$ (which might be equal to $r$) plays the role of a fixed reference probability distribution to be compared empirically with the probability distribution $r$.\\
The proof is again an application of the delta method in conjunction with the limit law from \Cref{thm:optimalsolutionlimit}. We again do not derive any substantially different distributional limit behaviour between the cases $r=s$ and $r\neq s$. This is in notably contrast to the non-regularized OT (see the discussion in \Cref{Introduction} and \Cref{Scaling}).

\begin{theorem} \label{thm:optimaltransportlimit}
Under the assumptions of \Cref{thm:optimalsolutionlimit}, as $n \to \infty$ it holds that
\begin{equation*}
\sqrt{n}\left\lbrace W_{p,\lambda,f} (\hat{r}_n,s)-W_{p,\lambda,f}(r,s) \right\rbrace  \overset{D}{\longrightarrow} \mathcal{N}_1 \left(0,\sigma^2_{p,\lambda,f}(r\vert s)\right) 
\end{equation*}
with variance 
\begin{equation}
\sigma^2_{p,\lambda,f}(r\vert s)=\gamma^T\, \Sigma_{p,\lambda,f}(r \vert s) \, \gamma\, ,
\end{equation}
where $\gamma$ is the gradient of the function $\pi \mapsto \langle c_p, \pi \rangle^{\frac{1}{p}}$ evaluated at the regularized transport plan $\pi_{p,\lambda,f}(r,s)$, and $\Sigma_{p,\lambda,f}(r \vert s)$ is the covariance matrix from \Cref{thm:optimalsolutionlimit}.
Standardizing by the square root of the empirical variance $\sigma^2_{p,\lambda,f}(\hat{r}_n\vert s)$ results in a standard normal limit distribution.
\end{theorem}

As a corollary, we immediately obtain limit distributions for the empirical entropy ROT plan and the empirical Sinkhorn divergence.

\begin{corollary}[Sinkhorn Transport and Sinkhorn Divergence]\label{cor:sinkhorn}
Consider \\the negative Boltzmann-Shannon entropy $f$ in \eqref{eq:entropy}. Then the statements in \Cref{thm:optimalsolutionlimit} and \Cref{thm:optimaltransportlimit} remain valid. Note, that the gradient inherent in the corresponding covariance matrix \eqref{eq:covarianceregtransport} is given by Example \ref{example:entropy}.
\end{corollary}

\begin{remark}[Entropy ROT Type Functionals]
From \Cref{thm:optimalsolutionlimit} we easily derive asymptotic distributions for any sufficiently smooth function of the ROT plan. Exemplarily, we consider the objective function in \eqref{eq:EntropyROT} denoted as $d(\pi_{p,\lambda,f}(r,s))$. A straightforward calculation shows that
\begin{equation*}
\nabla d(\pi_{p,\lambda,f}(r,s))= c_p +\lambda \log(\pi_{p,\lambda,f}(r,s))=(\alpha_{p,\lambda,f},{\beta_{p,\lambda,f}}_{\star})A_{\star}\, .
\end{equation*}
The second equality follows by primal-dual optimality relation between $\pi_{p,\lambda,f}$ and its optimal dual solutions $(\alpha_{p,\lambda,f},{\beta_{p,\lambda,f}}_{\star})$ \citep[Proposition 4.4]{peyre2017computational} with lower subscript star as we delete the last constraint in \eqref{eq:EntropyROT} (\Cref{rem:asymetric}). In conjunction with \Cref{example:entropy} and \Cref{thm:optimalsolutionlimit} we conclude
\begin{equation}\label{eq:objectivelimit}
\sqrt{n}\left\lbrace d(\pi_{p,\lambda,f}(\hat{r}_n,s))- d(\pi_{p,\lambda,f}(r,s)) \right\rbrace \overset{D}{\longrightarrow} \left\langle G,\alpha_{p,\lambda,f} \right\rangle\, ,
\end{equation}
where $G\sim \mathcal{N}_N \left(0,\Sigma(r)\right)$ \citep[Theorem 2.5]{bigot2017central}. Notably, if $r=s$ the limit law in \eqref{eq:objectivelimit} is non degenerate. This is not true anymore for the Sinkhorn loss \citep{genevay2017learning} defined by
\begin{equation*}
S_{\lambda}(r,s)\coloneqq d(\pi_{p,\lambda,f}(r,s)) - \frac{1}{2}\left(d(\pi_{p,\lambda,f}(r,r))-d(\pi_{p,\lambda,f}(s,s))\right)\, ,
\end{equation*}
as then $\nabla S_{\lambda}(r,r)=0$ (see also \cite{bigot2017central}). However, a second order expansion which is based on a perturbation analysis for the dual solutions provides a non degenerate asymptotic limit of $nS_{\lambda}(\hat{r}_n,r)$. This can be represented as a weighted sum of independent $\chi_1^2$ random variables. Exact computation is tedious but follows along the lines of \cite{luise2018differential} who also provide a perturbation analysis for the dual solutions. The weights of this sum are then given by the eigenvalues of the Hessian $\nabla^2 S_{\lambda}(r,r)$. From this it can be shown that the $m$ out of $n$ bootstrap is consistent when $m=o(n)$ \citep{shao1994bootstrap,rippl2016limit} which is an alternative to the bootstrap suggested in \cite{bigot2017central}.
\end{remark}

\subsection{Estimating Both Probability Distributions}\label{Estimationboth}
\Cref{thm:optimalsolutionlimit} and \Cref{thm:optimaltransportlimit} also hold true if we estimate both distributions $r$ and $s$ by their empirical versions $\hat{r}_n$ and $\hat{s}_m$ in \eqref{eq:empiricalmeasures} derived by two collections of $\mathcal{X}$-valued random variables $X_1,\ldots,X_n \overset{i.i.d.}{\sim} r$ and independently $Y_1,\ldots,Y_m \overset{i.i.d.}{\sim} s$, respectively. Note that we treat $r$ and $s$ asymmetrically (see \Cref{rem:asymetric}). Hence, the underlying multinomial process is based on the reduced multinomial vector $(r,s_{\star})$ rather than $(r,s)$. The scaling rate is given by $\sqrt{\nicefrac{nm}{n+m}}$ such that $n \wedge m \to +\infty$ and $\nicefrac{m}{n+m}\to \delta \in (0,1)$.\\
For instance, the limit distribution for the empirical ROT plan with parameters $\lambda>0$, $p\geq 1$ and proper regularizer $f$ then reads as 
\begin{equation*}
\sqrt{\frac{nm}{n+m}} \left\lbrace \pi_{p,\lambda,f}(\hat{r}_n,\hat{s}_m)-\pi_{p,\lambda,f}(r,s) \right\rbrace  \overset{D}{\longrightarrow} \mathcal{N}_{N^2}\left(0,\Sigma_{p,\lambda,f}(r, s)\right)\, .
\end{equation*}
The variance $\Sigma_{p,\lambda,f}(r, s)$ is different to $\Sigma_{p,\lambda,f}(r\vert s)$ from \Cref{thm:optimalsolutionlimit}. More precisely, given the covariance matrix of the reduced multinomial process
\begin{equation*}
\Sigma(\delta,r,s_{\star}) = \begin{pmatrix}
\delta\, \Sigma(r) & 0 \\ 
0 & (1-\delta)\,\Sigma(s_{\star})
\end{pmatrix} \, ,
\end{equation*}
we find that
\begin{equation}\label{eq:twosamplecovariance}
\Sigma_{p,\lambda,f}(r, s) = \nabla \phi_{p,\lambda,f}(r,s_{\star})\, \Sigma(\delta,r,s_{\star})\, \nabla \phi_{p,\lambda,f}(r,s_{\star})^T\, .
\end{equation}
Similar the limit distribution for the empirical ROT distance now reads as
\begin{equation*}
\sqrt{\frac{nm}{n+m}} \left\lbrace W_{p,\lambda,f}(\hat{r}_n,\hat{s}_m) - W_{p,\lambda,f}(r,s) \right\rbrace \overset{D}{\longrightarrow} \mathcal{N}_1 \left(0,\sigma^2_{p,\lambda,f}(r,s)\right) \, .
\end{equation*}
The variance $\sigma^2_{p,\lambda,f}(r,s)$ is again different to $\sigma^2_{p,\lambda,f}(r\vert s)$ from \Cref{thm:optimaltransportlimit} and given by $\sigma^2_{p,\lambda,f}(r,s) = \gamma^T\, \Sigma_{p,\lambda,f}(r, s)\, \gamma$, where we recall the definition of $\gamma$ from \Cref{thm:optimaltransportlimit}. 

\subsection{Comparison to (non-regularized) Optimal Transport} \label{Scaling}
The Sinkhorn divergence approximates the OT distance exponentially fast as $\lambda$ tends to zero. More precisely, there exists a constant $C$ depending on the support of $r$ and $s$ and the cost $c_p$ but independent of $\lambda$ such that $\left\vert W^p_{p,\lambda,f}(r,s)-W^p_p(r,s)\right\vert \leq C \exp(-\frac{1}{\lambda})$ \citep[Proposition 1]{luise2018differential} for $f$ the Boltzmann-Shannon entropy. For the finite setting considered here we obtain
\begin{equation}\label{eq:ROTbound}
\sup_{r,s \in \Delta_N}\left\vert W^p_{p,\lambda,f}(r,s)-W^p_p(r,s)\right\vert \leq C \exp\left(-\frac{1}{\lambda}\right)\, ,
\end{equation}
hence independent of the underlying (possibly unknown) distributions $r$ and $s$. By a similar argument as in \citet[Theorem 2.10]{bigot2017central} we decompose 
\begin{equation}
\begin{split}
&\sqrt{n}\left\lbrace W^p_{p,\lambda(n),f} (\hat{r}_n,s)-W^p_{p,\lambda(n),f}(r,s) \right\rbrace\nonumber\\ = &\sqrt{n}\left\lbrace W^p_{p,\lambda(n),f} (\hat{r}_n,s)-W^p_{p}(\hat{r}_n,s) \right\rbrace+\sqrt{n}\left\lbrace W^p_{p} (r,s)-W^p_{p,\lambda(n),f}(r,s) \right\rbrace\label{eq:lamrate}\\ &+ \sqrt{n}\left\lbrace W^p_{p} (\hat{r}_n,s)-W^p_{p}(r,s) \right\rbrace\, .\nonumber
\end{split}
\end{equation}
and obtain from \eqref{eq:lamrate} for $\lambda(n)=o\left(\nicefrac{1}{\log(\sqrt{n})}\right)$ that the first two terms converge to zero while the third term asymptotically follows the limit law for the (non-regularized) transport distance \citep{sommerfeld2018inference}. Hence, if $\lambda(n)$ converges faster to zero as $\nicefrac{1}{\log(\sqrt{n})}$ the limit law of the empirical Sinkhorn divergence in \eqref{eq:lamrate} is the same as for its non-regularized counterpart. Note that this holds for $W^p_{p,\lambda,f}$ and in fact fails for $W_{p,\lambda,f}$ in the case $r=s$. Further, this provides a notable difference between the logarithmic rate for $\lambda(n)$ obtained here and the rate $\lambda(n)=o\left(\nicefrac{1}{\sqrt{n}}\right)$ by \citet{bigot2017central} considering the optimal value of \eqref{eq:EntropyROT}, as the true limit for $\lambda=0$ is obtained for a much larger range of $\lambda(n)$-sequences as for the optimal value of \eqref{eq:EntropyROT}. This indicates that computation of the Sinkhorn divergence as defined in \eqref{eq:sinkhorn} provides a numerically more stable approximation to the OT for small $\lambda$ as for the optimal value in \eqref{eq:EntropyROT}, supporting our numerical findings in \Cref{Simulations}.

\section{Bootstrap} \label{Bootstrap}

The (non-regularized) OT distance on finite spaces defines a functional that is only \textit{directionally} Hadamard differentiable, i.e. has a \textit{non-linear} derivative with respect to $r$ and $s$. Hence, the (naive) $n$ out of $n$ bootstrap method to approximate the distributional limits for the (non-regularized) empirical OT distance fails \citep{sommerfeld2018inference}. However, our arguments underlying the proof of \Cref{thm:optimalsolutionlimit} for ROT are based on the usual delta method for \textit{linear} derivatives. As a byproduct we obtain that for the ROT plan and for the ROT distance the $n$ out of $n$ bootstrap is consistent. More precisely, conditionally on the data $X_1,\ldots,X_n$, the law of the multinomial empirical bootstrap process $\sqrt{n}\{\hat{r}^*_n - \hat{r}_n\}$ is an asymptotically consistent estimator of the law for the multinomial empirical process $\sqrt{n}\{ \hat{r}_n - r\}$ \citep[Theorem 3.6.1]{van1996weak}. Here, $\hat{r}^*_n=\frac{1}{n}\sum_{i=1}^n \delta_{X^*_i}$ is the empirical bootstrap estimator for $\hat{r}_n$ derived by a sample $X^*_1,\ldots,X^*_n\overset{i.i.d.}{\sim} \hat{r}_n$. Such conditional weak convergence can be formulated in terms of the bounded Lipschitz metric, that is
\begin{equation}\label{eq:bootstrapmarginalone}
\sup_{h \in \text{BL}_1\left(\mathbb{R}^N\right)} \left\vert \mathbb{E}[h(\sqrt{n}\{\hat{r}^*_n - \hat{r}_n\}) \vert  X_1,\ldots,X_n] - \mathbb{E}[h(\sqrt{n}\{ \hat{r}_n - r\})] \right\vert
\end{equation}
converges to zero in probability, where
\begin{equation*}
\text{BL}_1\left(\mathbb{R}^{N}\right)\coloneqq \left\lbrace f\colon \mathbb{R}^{N} \to \mathbb{R} \, \bigr\vert \sup_{x \in \mathbb{R}^{N}}\lvert f(x)\rvert \leq 1, \, \lvert f(x_1)-f(x_2)\rvert \leq \Vert x_1-x_2 \Vert \right\rbrace
\end{equation*}
is the set of all bounded functions with Lipschitz constant at most one.
Combined with the consistency of the delta method for the bootstrap \citep[Theorem 3.9.11]{van1996weak}, this proves the consistency for the empirical bootstrap ROT plan. The statements again hold true for $r=s$ and $r\neq s$.

\begin{theorem}\label{thm:bootstrap}
Under the assumptions of \Cref{thm:optimalsolutionlimit} the (naive) bootstrap for the regularized optimal transport plan is consistent, that is
\begin{equation*}
\begin{split}
\sup_{h \in \text{BL}_1\left(\mathbb{R}^{N^2}\right)} \big\vert &\mathbb{E}[h(\sqrt{n}\left\lbrace \pi_{p,\lambda,f}(\hat{r}^*_n,s)-\pi_{p,\lambda,f}(\hat{r}_n,s) \right\rbrace  \vert  X_1,\ldots,X_n] \\
&- \mathbb{E}[h(\sqrt{n}\left\lbrace \pi_{p,\lambda,f}(\hat{r}_n,s)-\pi_{p,\lambda,f}(r,s) \right\rbrace )] \big\vert \overset{\mathbb{P}}{\longrightarrow} 0\, .
\end{split}
\end{equation*}
This holds as well for the regularized optimal transport distance
\begin{equation*}
\begin{split}
\sup_{h \in \text{BL}_1(\mathbb{R})} \big\vert &\mathbb{E}[h(\sqrt{n}\left\lbrace W_{p,\lambda,f}(\hat{r}^*_n,s)-W_{p,\lambda,f}(\hat{r}_n,s) \right\rbrace  \vert  X_1,\ldots,X_n] \\
&- \mathbb{E}[h(\sqrt{n}\left\lbrace W_{p,\lambda,f}(\hat{r}_n,s)-W_{p,\lambda,f}(r,s) \right\rbrace )] \big\vert \overset{\mathbb{P}}{\longrightarrow} 0\, .
\end{split}
\end{equation*}
\end{theorem}

Analogously, the bootstrap consistency is valid if $Y^*_1,\ldots,Y^*_m \overset{i.i.d.}{\sim} \hat{s}^*_m$ independently to $X^*_1,\ldots,X^*_n \overset{i.i.d.}{\sim} \hat{r}^*_n$ (see \Cref{Distributions}).

\section{Simulations}\label{Simulations}

We illustrate our distributional limit results in Monte Carlo simulations. Exemplarily, we investigate the speed of convergence for the empirical Sinkhorn divergence $(p=1)$ to its limit distribution (Corollary \ref{cor:sinkhorn}) in the one-sample case in both settings $r=s$ and $r \neq s$. Moreover, we illustrate the accuracy of approximation by the (naive) $n$ out of $n$ bootstrap (Theorem \ref{thm:bootstrap}). As the Sinkhorn divergence approximates the (non-regularized) OT distance, we also compare for small regularization parameters the finite sample distribution of the Sinkhorn divergence with the limit laws for the (non-regularized) optimal transport distance (OT distance) in \eqref{eq:finitetransportlimit}.\\
All simulations were performed using \textsf{R} \citep{R}. The Sinkhorn divergences are calculated with the \textsf{R}-package \textit{Barycenter} \citep{package:barycenter}.

\begin{remark}[Computation of (empirical) variances]
For $p=1$ it holds that $\sigma^2_{1,\lambda,f}(r \vert s)=c_1^T\, \Sigma_{1,\lambda,f}(r\vert s)\, c_1$. According to \Cref{example:entropy} and \Cref{cor:sinkhorn}, the computation of the variance involves the computation of 
\begin{equation*}
\begin{split}
\Sigma_{1,\lambda,f}(r\vert s) &= \nabla_r \phi_{1,\lambda,f}(r,s_{\star}) \, \Sigma(r)\, [\nabla_r \phi_{1,\lambda,f}(r,s_{\star})]^T \\
&= D\,A_{\star}^T [A_{\star}\,D\, A_{\star}^T]^{-1}_{[1:N]}\, \Sigma(r) \, [A_{\star}\,D\, A_{\star}^T]^{-1}_{[N:1]}\, A_{\star}\,D \, ,
\end{split}
\end{equation*}
where the subscript notation $[1:N]$ $([N:1])$ denotes the first $N$ columns (rows) of the corresponding matrix. Recall that $D$ is equal to a diagonal matrix with diagonal given by the entropy ROT plan $\pi_{1,\lambda,f}(r,s)$ and $A_{\star}$ is the coefficient matrix in \eqref{eq:OT} reduced by its last row. Besides calculating the entropy ROT plan, the computation of the variance faces matrix inversion. However, the matrix $A_{\star}\,D\,A_{\star}^T$ is symmetric and positive definite by full rank of $A_{\star}$. Moreover, it posses a block structure given by
\begin{equation*}
A_{\star}\,D\,A_{\star}^T = \begin{bmatrix}
R & \Pi_{1,\lambda,f}\\ \Pi_{1,\lambda,f}^T & S_{\star}
\end{bmatrix}\, ,
\end{equation*}
where $\Pi_{1,\lambda,f}$ denotes the matrix version of the entropy ROT plan reduced by its last column, and we set $R\coloneqq \text{diag}(r)$ and $S_{\star} \coloneqq \text{diag}(s_{\star})$. Hence, we can apply a block wise inversion and obtain
\begin{equation*}
[A_{\star}\,D\,A_{\star}^T]_{[\,,1:N]}^{-1} = \begin{bmatrix}
[R-\Pi_{1,\lambda,f}\, S_{\star}^{-1}\, \Pi_{1,\lambda,f}^T]^{-1} \\
-S_{\star}^{-1}\, \Pi_{1,\lambda,f}^T[R-\Pi_{1,\lambda,f}\,S_{\star}^{-1}\,\Pi_{1,\lambda,f}^T]^{-1} 
\end{bmatrix}\, .
\end{equation*}
\end{remark}

We consider the finite metric space $\mathcal{X}$ to be an equidistant two-dimen\-sional $L \times L$ grid on $[0,1]\times[0,1]$ and the cost $c \in \R^{L^4}$ consisting of the euclidean distance $(p=1)$ between the pixels on the grid. Note that different grid sizes for fixed regularization parameter $\lambda$ correspond to different amounts of regularization. As recommended by \citet{cuturi2013sinkhorn}, we let the amount of regularization depend on the median distance $q_{50}(c)$ between the pixels on the grid. More precisely, we define the regularization parameter by
\begin{equation}\label{eq:lambdaregularization}
\lambda \coloneqq \lambda_0 \, q_{50}(c)\, ,
\end{equation}
where $\lambda_0>0$ is a parameter that we vary for different simulations. The probability distributions $r,s \in \Delta_{L^2}$ on $\mathcal{X}$ are generated as two independent realizations of a Dirichlet random variable $\text{Dir}(\boldsymbol{\alpha})$ with concentration parameter $\boldsymbol{\alpha}=(\alpha, \ldots, \alpha) \in \R^{L\times L}$. The choice $\alpha=1$ corresponds to a uniform distribution on the probability simplex $\Delta_{L^2}$.

\subsection{Speed of Convergence}

We first generate for grid size $L=10$ probability distributions $r$ on $\mathcal{X}$ as independent realizations of a $\text{Dir}(\boldsymbol{1})$ random variable. Given such a distribution, we fix the amount of regularization to $\lambda=2\,q_{50}(c)$, that is $\lambda_0=2$ in \eqref{eq:lambdaregularization}. We then sample $n=25$ observations according to this probability distribution $r$ and compute 
\begin{equation*}
\sqrt{\frac{n}{\sigma^2_{1,\lambda,f}(\hat{r}_n \mid r)}}\left\lbrace W_{1,\lambda,f} (\hat{r}_n,r)-W_{1,\lambda,f}(r,r) \right\rbrace \, ,
\end{equation*}
referred to as a Sinkhorn sample. This is repeated $20.000$ times and simulates the scenario when the data generating probability distribution $r$ coincides with the probability distribution to be compared. Similarly, we consider the same set up in the case $r\neq s$ when we simulate independently a second distribution $s\sim \text{Dir}(\boldsymbol{1})$. The finite sample distributions are then compared to their theoretical limit distributions which by Theorem \ref{thm:optimaltransportlimit} are standard normal distributions.

\begin{figure}
  \centering
  \setlength\tabcolsep{8pt}
  \begin{tabular}{cc}
    \includegraphics[width=0.45\textwidth]{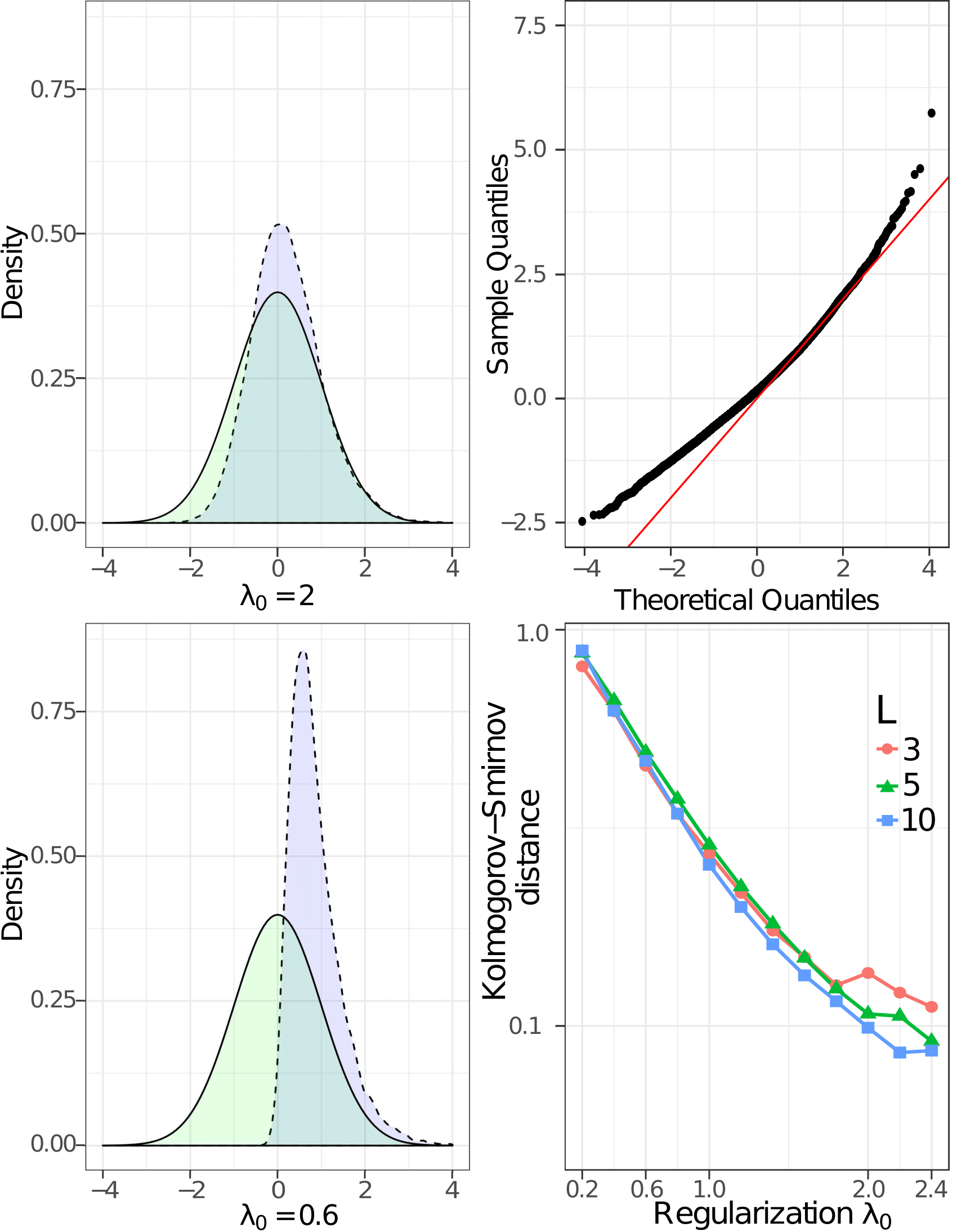} & \includegraphics[width=0.45\textwidth]{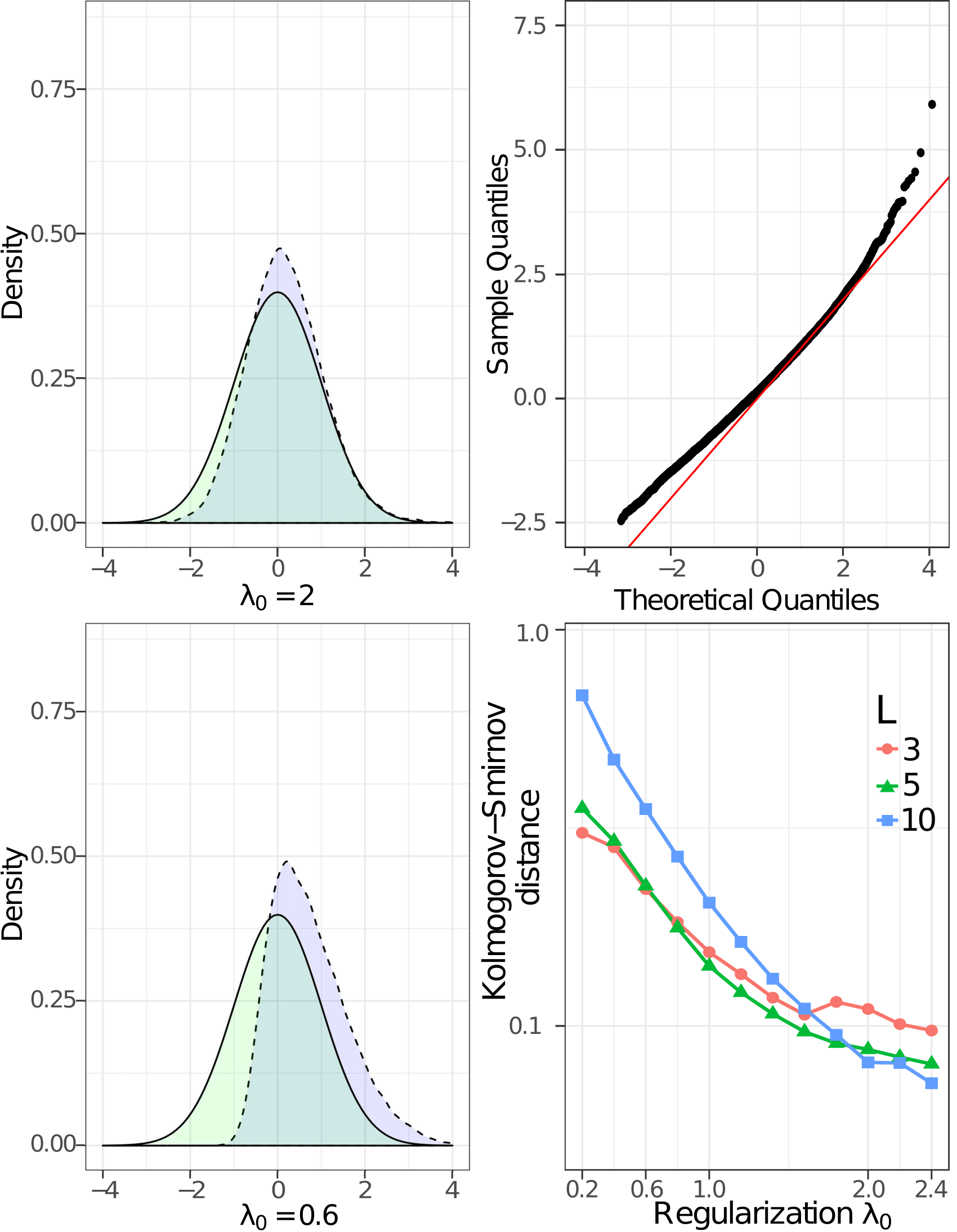} \\
     (a) $r=s$ &  (b) $r\neq s$ 
  \end{tabular}
  \caption{\label{fig:densityqqplot}\textbf{(a) Finite sample accuracy of the limit law in the one sample case for $\mathbf{r=s}$.} First row: Finite sample density (dashed line) of the empirical Sinkhorn divergence for $n=25$ on a regular grid of size $L=10$ with regularization parameter $\lambda_0=2$ compared to its limit law (standard Gaussian, solid line). The finite sample density has been estimated with a kernel density estimator with Gaussian kernel and Silverman's rule to select bandwidth. On the right the corresponding Q-Q-plot, where perfect fit is indicated by the red solid line.
Second row: L.h.s. same setting as above, $\lambda_0=0.6$. R.h.s. Finite sample accuracy in dependence on $\lambda_0$: The Kolmogorov-Smirnov distance on a logarithmic scale averaged over five realization of a $\text{Dir}(\boldsymbol{1})$ distribution between the finite sample distribution $(n=25)$ of the empirical Sinkhorn divergence and the standard normal distribution as a function of the regularization parameter $\lambda_0$.
 \newline \textbf{(b) Finite sample accuracy of the limit law in the one sample case for $\mathbf{r\neq s}$.} Same scenario as in (a) whereas here the probability distribution $r$ to be sampled is not equal to the fixed reference probability distribution $s$.}
\end{figure}

\begin{figure}[ht]
  \centering
  \begin{tabular}{cc}
    \includegraphics[width=0.45\textwidth]{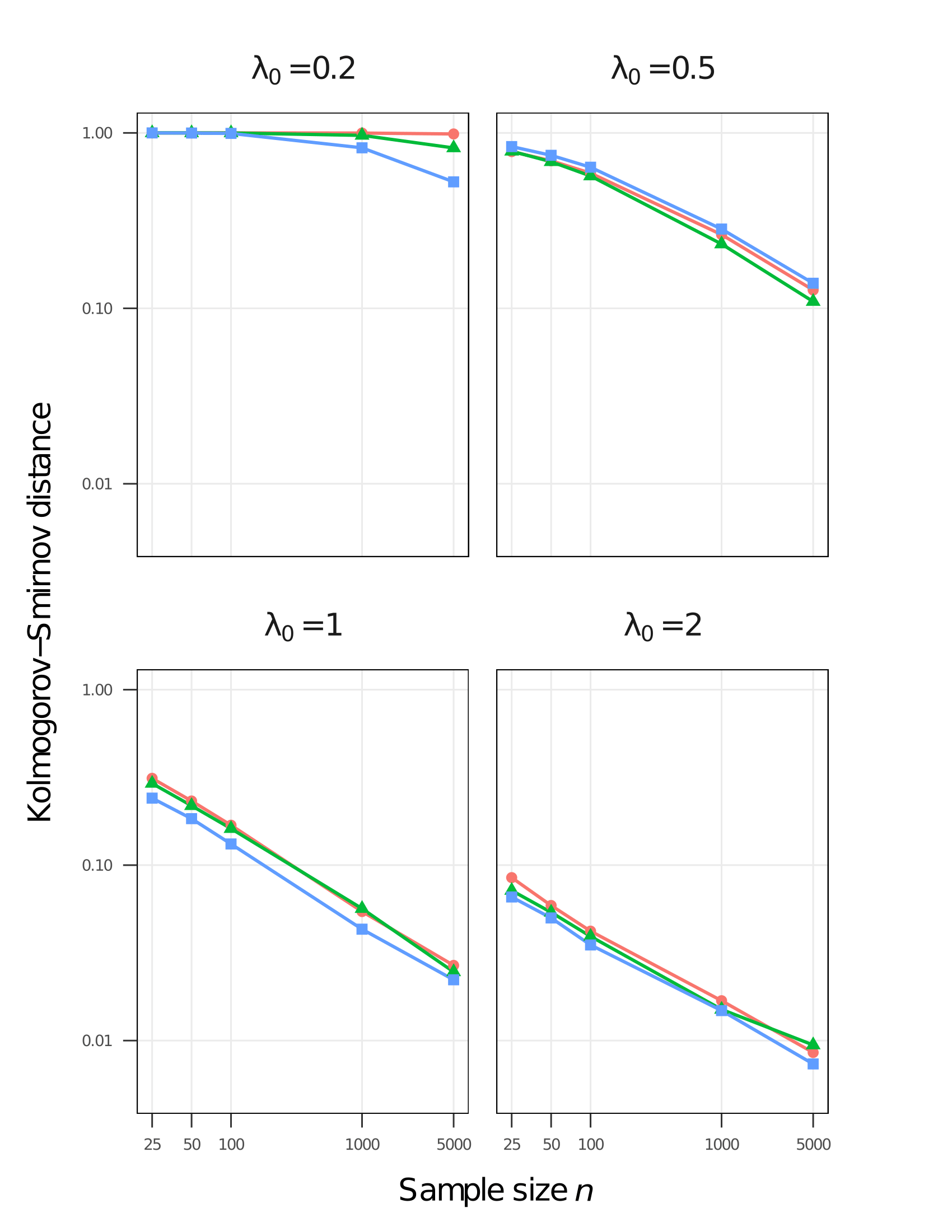} & \includegraphics[width=0.45\textwidth]{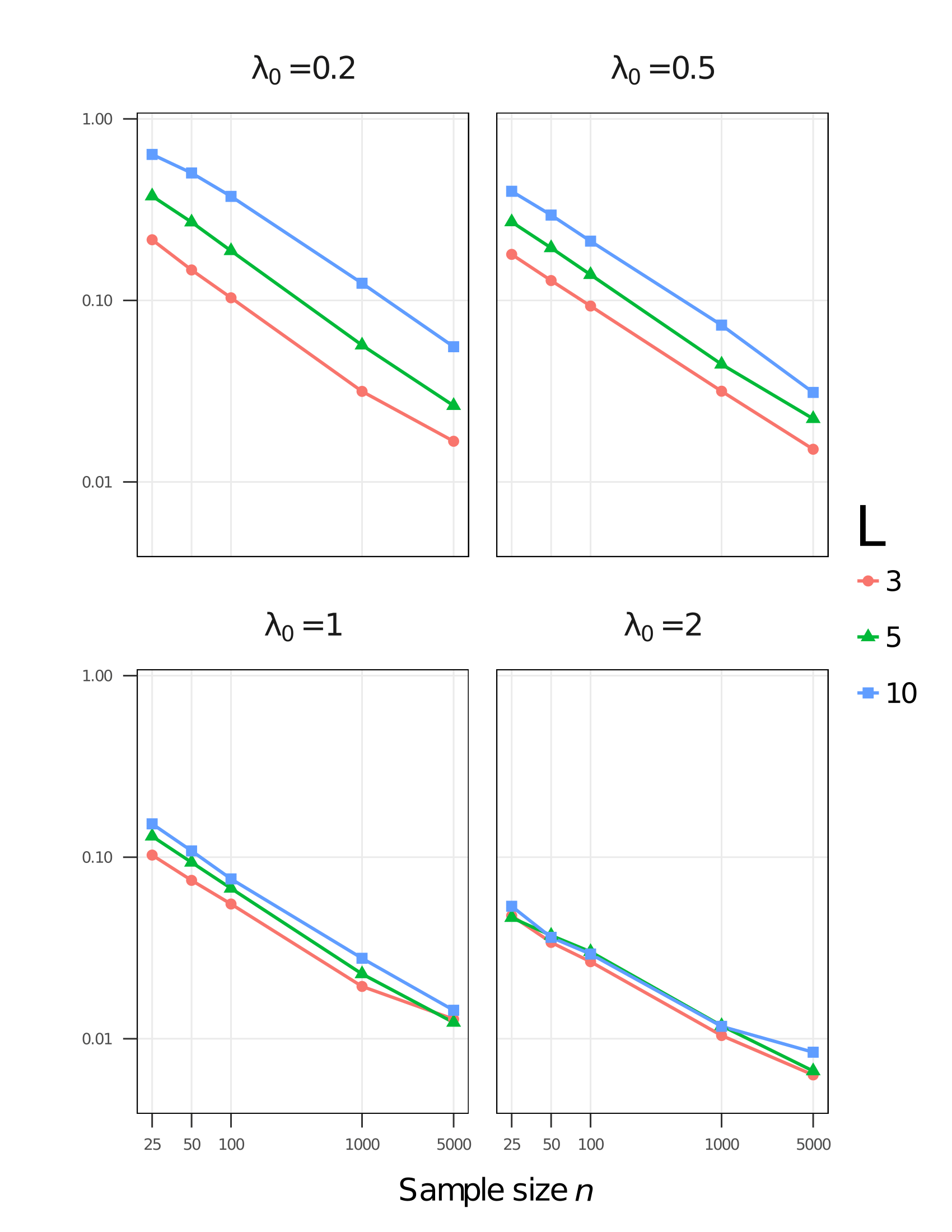} \\
     (a) $r=s$ &  (b) $r\neq s$ 
  \end{tabular}
\caption{\textbf{(a) Kolmogorov-Smirnov distance in the one-sample case for $\mathbf{r=s}$} The Kolmogorov-Smirnov distance between the finite Sinkhorn divergence sample distribution and its theoretical normal distribution for $r=s$ (left) and $r\neq s$ (right) as a function of the sample size $n\in \{25,\,50,\,100,\,1000,\,5000\}$ for different grid sizes $L\times L$ and different regularization parameters $\lambda_0$. The distances are averaged over five (pairs) of independent realizations of a $\text{Dir}(\boldsymbol{1})$ distribution. The axes are given on a logarithmic scale. \newline
\textbf{(b) Kolmogorov-Smirnov distance in the one-sample case for $\mathbf{r\neq s}$} Same scenario as in (a) whereas here the probability distribution $r$ to be sampled is not equal to the fixed reference probability distribution $s$.}
\label{fig:KSdistance}
\end{figure}

We demonstrate the results by kernel density estimators and corresponding Q-Q-plots in the first row of Figure \ref{fig:densityqqplot} (a) and (b). We observe that the finite sample distribution is already well approximated for small sample size $(n=25)$ by the theoretical limit distribution.
However, the amount of entropy regularization $(\lambda_0=2)$ added to OT is rather large. We find that for sample size $n=25$ the smaller the regularization $\lambda_0$ the worse the approximation by the theoretical Gaussian limit law. This is depicted in the second row of Figure \ref{fig:densityqqplot} (a) and (b) where we analyse for small regularization parameters the Kolmogorov-Smirnov distance (maximum absolute difference between empirical cdf and cdf of the limit law). Note, that the theoretical limit law approximates the finite sample distribution slightly better in the case $r\neq s$.\\
We additionally investigate the speed of convergence with respect to the Kolomogorov-Smirnov distance of the finite sample distribution in the small regularization regime when the sample size is large $(n\gg 25)$. As illustrated in Figure \ref{fig:KSdistance} we observe in our simulations good approximation results for rather large regularization whereas for small regularization the approximation accuracy decreases. This can only be compensated if the sample size $n$ severely increases to several thousands. As a benchmark, for $\lambda_0=0.2$ we already require $n=5000$ samples to observe an accurate approximation by the limit distribution from Theorem \ref{thm:optimaltransportlimit}.

\begin{figure}[ht]
  \centering
  \setlength\tabcolsep{8pt}
  \begin{tabular}{cc}
    \includegraphics[width=0.45\textwidth]{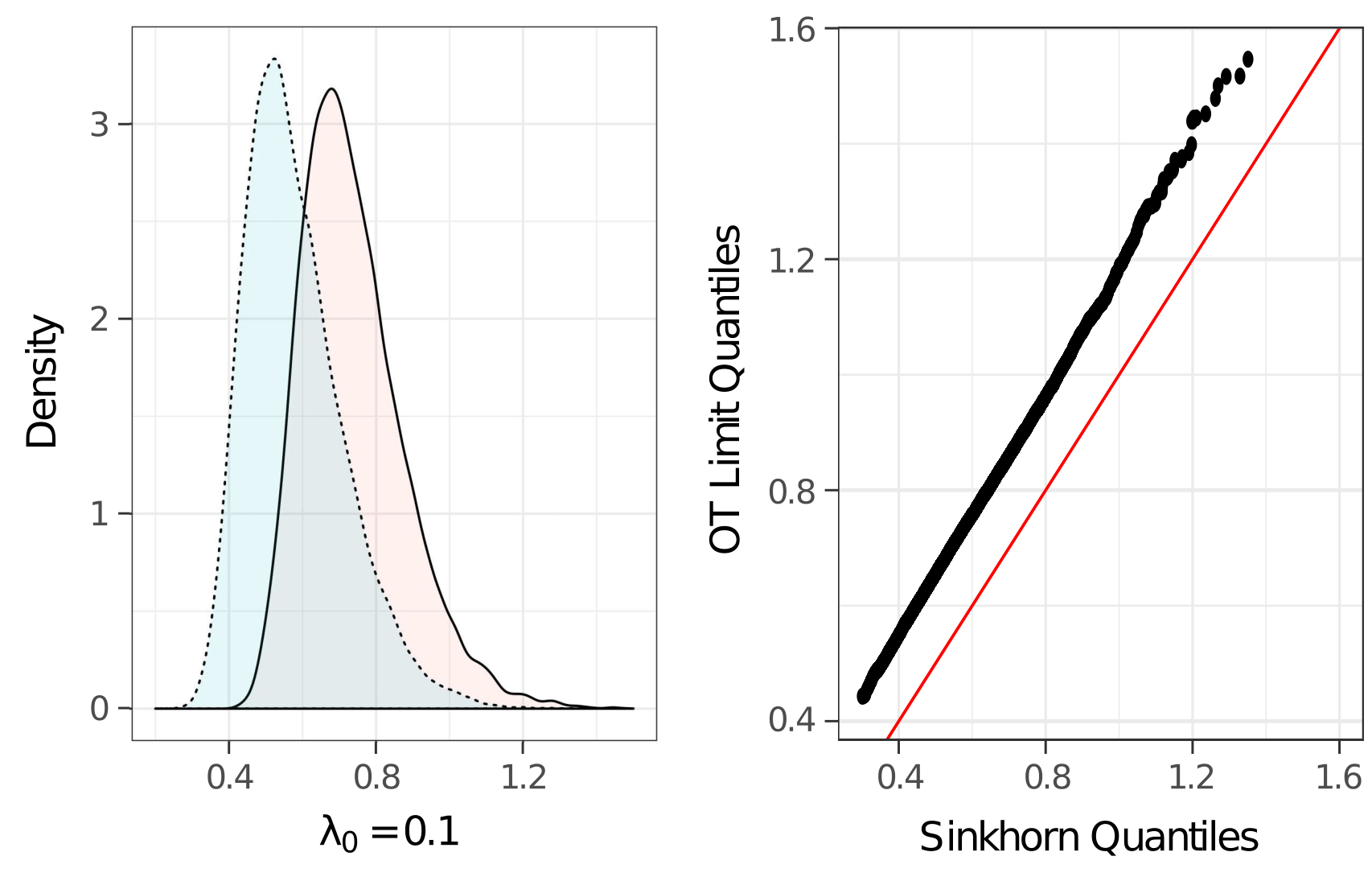} & \includegraphics[width=0.45\textwidth]{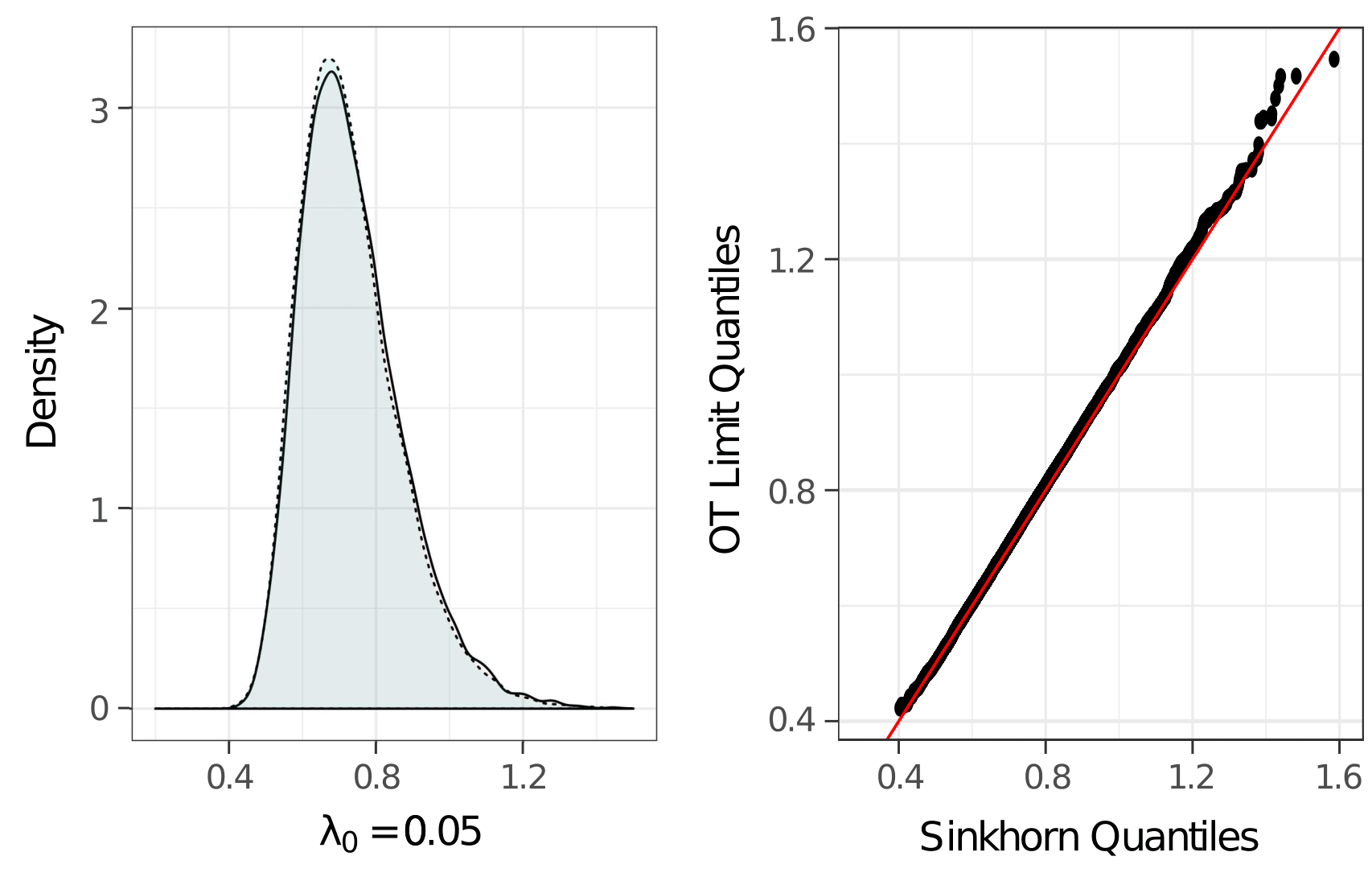}
  \end{tabular}
\caption{\textbf{Comparison to OT limit law in the one-sample case for $\mathbf{r=s}$.} Comparison of the finite sample distribution $(n=25)$ of the empirical Sinkhorn divergence on an equidistant grid of size $L=10$ for regularization parameter $\lambda_0=0.1$ (left two figures) and $\lambda_0=0.05$ (right two figures) to the OT limit law in \eqref{eq:finitetransportlimit} \citep[Theorem 1]{sommerfeld2018inference}. Kernel density estimator for the Sinkhorn sample (dotted line) and the OT sample (solid line) with corresponding Q-Q-plot on the right. The OT distance limit distribution is approximated by a sample of size $20.000$ implemented in the \textsf{R}-package \textit{otinference} \citep{package:otinference}.}
\label{fig:OTcomparison}
\end{figure}

Motivated by the inaccurate approximation for small sample size in the small regularization regime and the fact that the Sinkhorn divergence converges to the OT distance as $\lambda_0 \searrow 0$, we compare the finite sample distribution to the (non-regularized) OT limit law in \citet{sommerfeld2018inference} for the case $r=s$. From Figure \ref{fig:OTcomparison} we see that the finite sample distribution for $\lambda_0=0.1$ and $\lambda_0=0.05$ and small sample size $(n=25)$ is approximated by the OT limit law in \eqref{eq:finitetransportlimit}. 

In summary, the finite sample distribution converges to its theoretical limit law. The accuracy of the approximation is driven by the regularization parameter $\lambda_0$. For large regularization added to OT, the limit law serves as a good approximation to the finite sample distribution, already for small sample sizes and independent of the size of the ground space. As $\lambda_0$ decreases, accuracy of approximation decreases which is consistent with our theoretical findings in \Cref{Scaling}.

\subsection{Simulating the Bootstrap}

\begin{figure}
  \centering
  \begin{tabular}{cc}
    \includegraphics[width=0.45\textwidth]{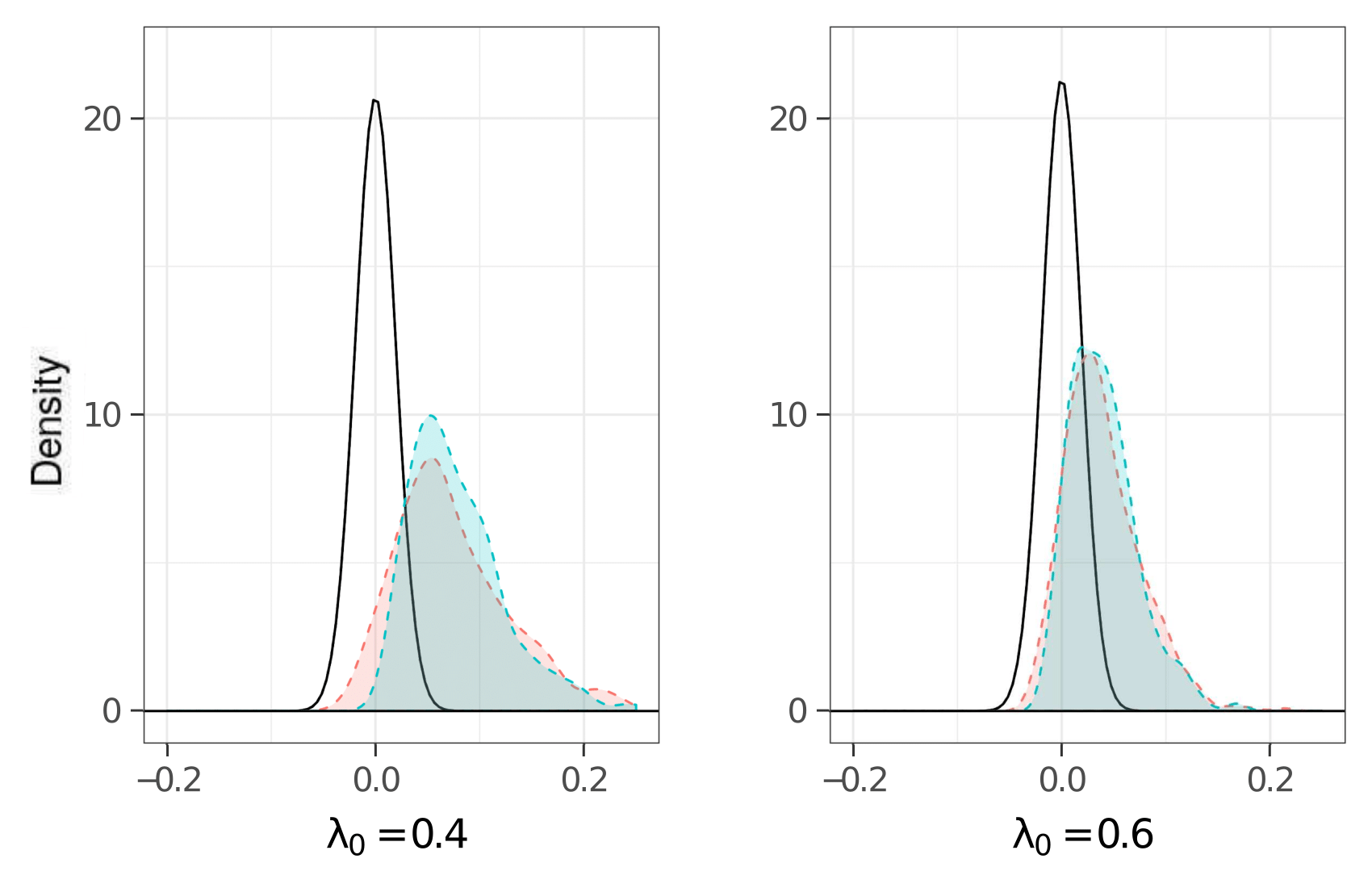} & \includegraphics[width=0.45\textwidth]{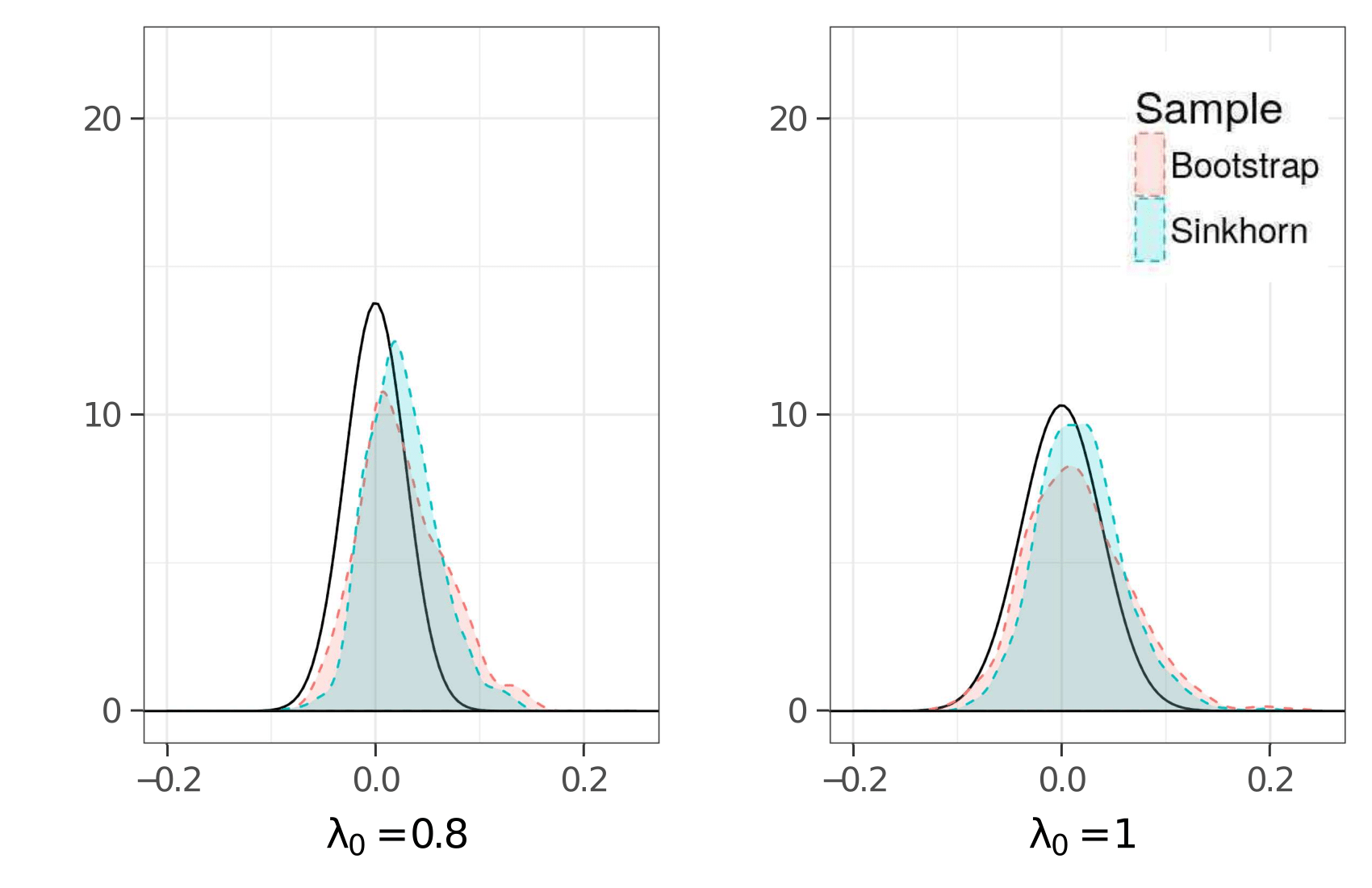} 
  \end{tabular}
\caption[Bootstrap in the one-sample case for $r=s$.]{\textbf{Bootstrap in the one-sample case for \boldmath$r=s$.} Illustration of the (naive) $n$ out of $n$ plug-in bootstrap approximation $(n=100)$ for different regularization parameters $\lambda_0$ on a grid of size $L=10$. The density in blue (resp. red) is obtained by a kernel density estimator (Gaussian kernel with bandwidth given by Silverman's rule) of a Sinkhorn divergence sample ($20.000$ realizations) \eqref{eq:sinkhornwassersteinsample} (resp. bootstrap sample \eqref{eq:sinkhornwassersteinbootstrapsample}, $B=500$ replications). The density of the corresponding Gaussian limit is depicted in black.}
\label{fig:bootstrap}
\end{figure}

Additionally, we simulate the (naive) $n$ out of $n$ plug-in bootstrap approximation from Section \ref{Bootstrap} in the one-sample case for $r=s$. For a grid with $L=10$ and cost given by the euclidean distance as before, we simulate $r \sim \text{Dir}(\boldsymbol{1})$ and generate $20.000$ realizations
\begin{equation}\label{eq:sinkhornwassersteinsample}
\sqrt{n}\left\lbrace W_{1,\lambda,f} (\hat{r}_n,r)-W_{1,\lambda,f}(r,r) \right\rbrace \, ,
\end{equation}
where we set the sample size $n=100$ and as before vary $\lambda_0$ (see \eqref{eq:lambdaregularization}). Moreover, for fixed empirical distribution $\hat{r}_n$ and each $\lambda_0$ we generate $B=500$ bootstrap replications 
\begin{equation}\label{eq:sinkhornwassersteinbootstrapsample}
\sqrt{n}\left\lbrace W_{1,\lambda,f} (\hat{r}_n^*,r)-W_{1,\lambda,f}(\hat{r}_n,r) \right\rbrace 
\end{equation}
by drawing independently with replacement $n=100$ times according to $\hat{r}_n$. We then compare the finite sample distributions again by kernel density estimators. The results are depicted in \Cref{fig:bootstrap}.

In fact, the finite bootstrap sample distribution is a good approximation of the finite sample distribution. However, as before, the speed of convergence to the corresponding limit distribution is driven by the amount of regularization. Similar results hold for the two-sample case (not displayed).

\section{Reducing Computational Complexity by Resampling}\label{Subsampling}

\citet{dvurechensky2018computational} analyse algorithms to compute the entropy ROT that yield $\epsilon$-approximates to the OT distance, i.e. $W^p_{p,\lambda,f}(r,s)\leq W^p_{p}(r,s) + \epsilon$. These methods are usually based on matrix scaling of the underlying $N\times N$ distance matrix in order to find the entropy ROT plan $\pi_{p,\lambda,f}(r,s)$ in \eqref{eq:ROTsolution}. With increasing number of support points $N$ of the underlying distributions, matrix scaling becomes computational infeasible mainly because of the high memory demand to store and scale the distance matrix.\\
In order to maintain computational feasibility, we modify an idea introduced by \cite{sommerfeld2018optimal}, i.e. we use a resampling scheme of the underlying probability distribution. The subsequent data example in \Cref{Colocalization} requires to deal with probability distributions with support on up to $N=300,000$ points (images represented by normalized gray scale pixel intensities). This requires storing a distance matrix with entries $300,000^2= 9\cdot 10^{12}$ which is far beyond the storage capacity of any standard laptop. In order to stay within the scope of computational feasibility, we resample from these probability distributions $n=50,000$ data points which reduces the required storage to $50,000^2=2.5\cdot 10^{9}$ entries. Note that while the number of support points of the resampled probability distributions decreases just by a sixth, the memory demand for the corresponding distance matrix is reduced by three orders of magnitude. This keeps matrix scaling feasible.

\section{Colocalization Analysis}\label{Colocalization}

Complex protein interaction networks are at the core of almost all cellular processes. To gain insight into these networks a valuable tool is colocalization analysis of images generated by fluorescence microscopy aiming to quantify the spatial proximity of proteins \citep{adler_quantifying_2010, zinchuk_vadim_quantitative_2014,wang2017spatially}. With the advent of super-resolution microscopy (nanoscopy) nowadays protein structures at a size of a single protein can be discerned. This challenges conventional colocalization methods, e.g. based on pixel intensity correlation as there is no spatial overlap due to blurring anymore \citep{xu_resolution_2016}.

\begin{figure}
  \centering
  \begin{tabular}{@{}c@{}}
    \includegraphics[width=0.75\textwidth]{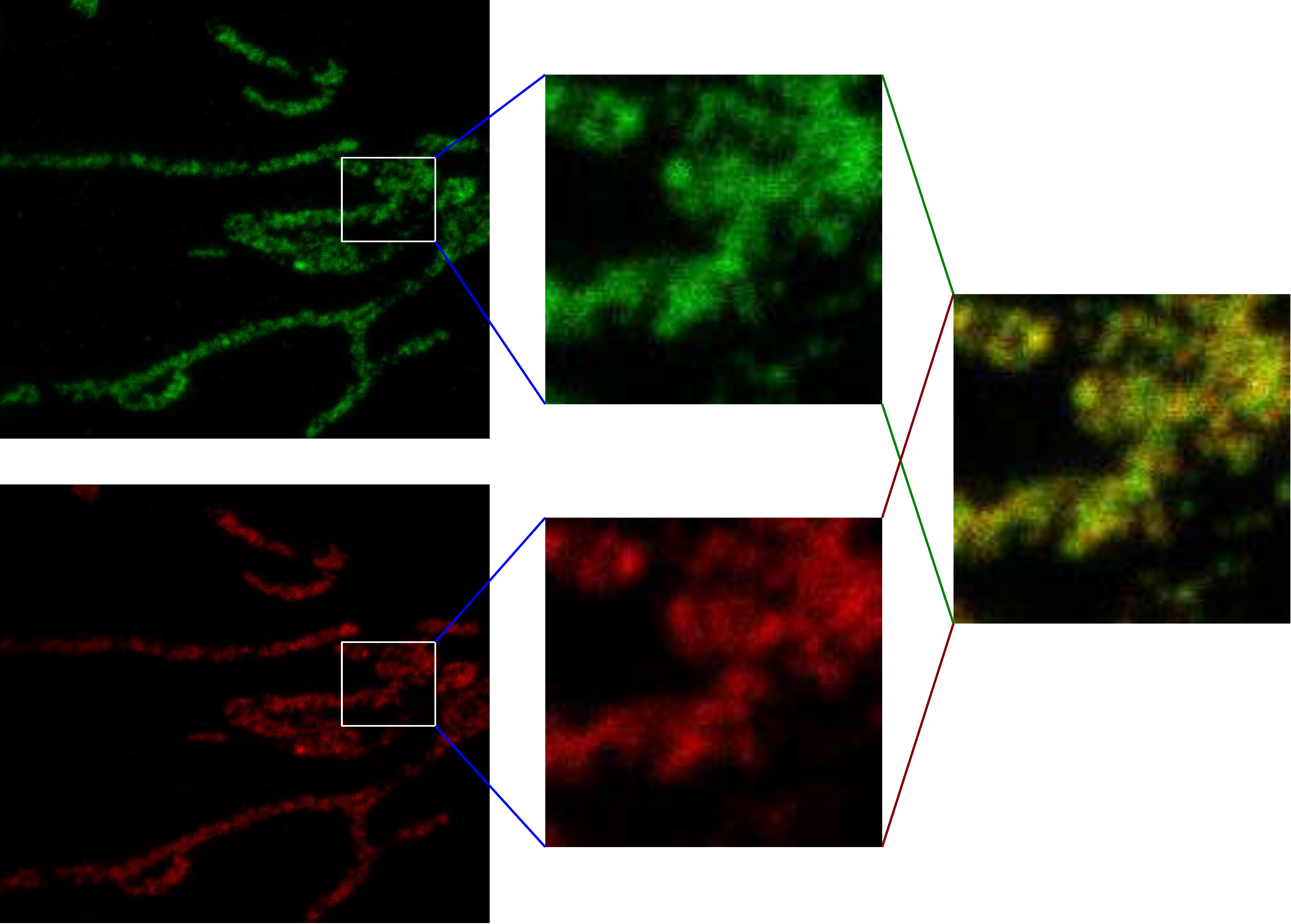} \\[\abovecaptionskip]
    \small  \textbf{(a) Setting 1:} Double staining of the protein \textit{ATP Synthase}.
  \end{tabular}

  \vspace{\floatsep}

  \begin{tabular}{@{}c@{}}
    \includegraphics[width=0.75\textwidth]{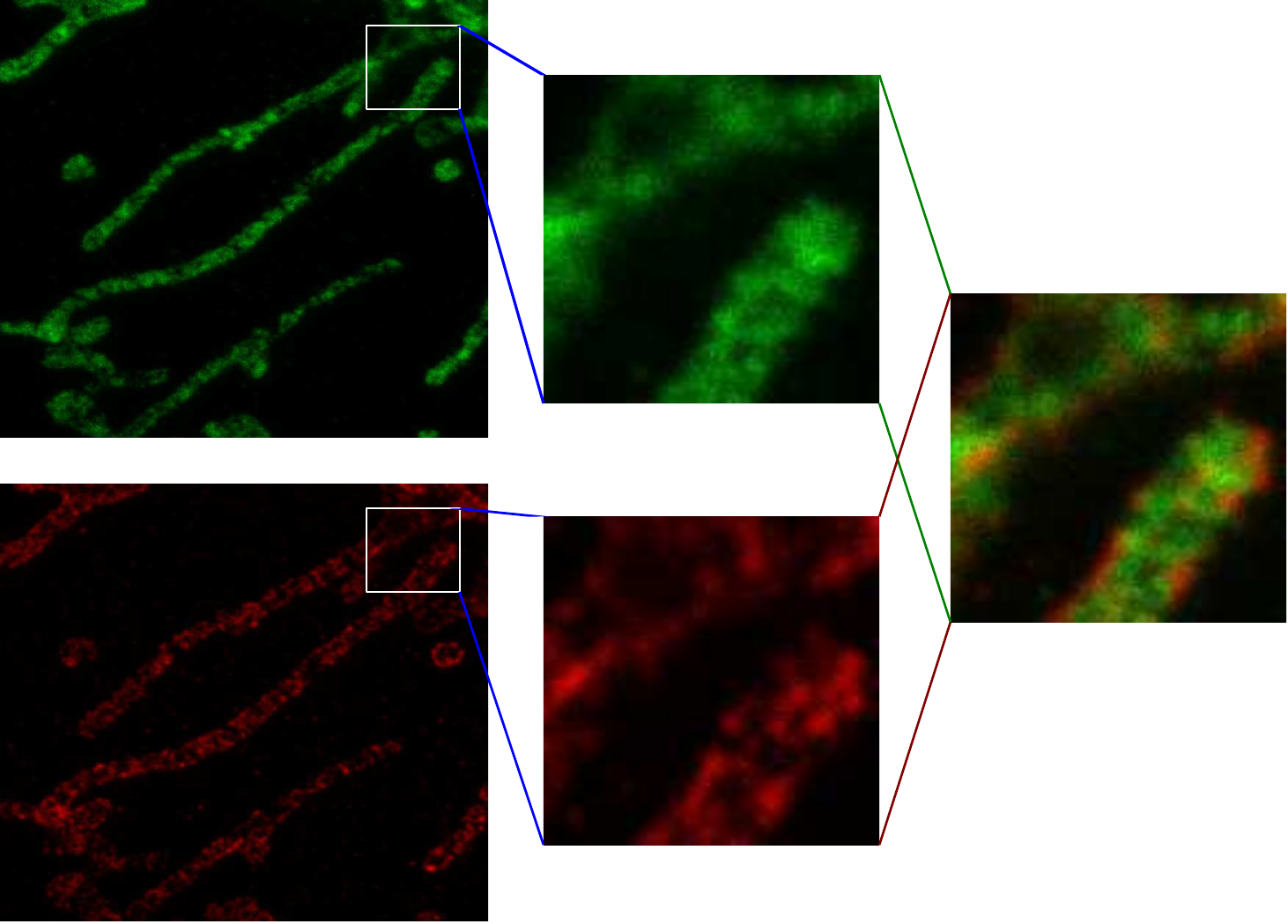} \\[\abovecaptionskip]
    \small  \textbf{(b) Setting 2:} Staining of the proteins \textit{ATP Synthase} (green) and \textit{MIC60} (red).
  \end{tabular}

  \caption{\label{fig:STED}\textbf{STED images for colocalization analysis.} Exemplary STED images of the two colocalization scenarios. Left: Images of a green and a red channel. Image size $666\times 666$ pixels, pixel size $=15$nm. Middle: Zoom ins ($128\times 128$ pixels). Right: Overlay of zoomed in images.}
\end{figure}

In \Cref{fig:STED} exemplary 2-colour-STED images (recorded at the Jakobs' lab, Department of NanoBiophotonics, Max-Planck Institute for Biophysical Chemistry, G\"ottingen) are displayed for illustrative purposes. These are generated by inserting two different fluorophore markers which emit photons at different wavelength (two colours) and reading them out, after excitation, with a stimulated emission depletion (STED) laser beam \citep{hell_far-field_2007}. \Cref{fig:STED} (a) shows 2-colour-STED images which were generated by attaching the two markers to the protein \textit{ATP Synthase}. \Cref{fig:STED} (b) displays the second case in which the markers are attached to two different proteins \textit{ATP Synthase} and \textit{MIC60}.\\
We aim to quantify the spatial proximity of \textit{ATP Synthase} and \textit{MIC60} (\Cref{fig:STED} (b)) and to set this in relation with the highest empirically possible colocalization represented by the double staining of the protein \textit{ATP Synthase} (\Cref{fig:STED} (a)). The overlay of the channels from setting 2 in \Cref{fig:STED} (b) already indicates that \textit{ATP Synthase} and \textit{MIC60} are located in different regions as there are only small areas which are yellow. In contrast, and as expected for the highest empirically possible colocalization, the yellow areas in the overlay of the two channels from the double staining in \Cref{fig:STED} (a) are more pronounced. 
In this section, we illustrate that the ROT plan \eqref{eq:ROTsolution} provides a useful tool to measure colocalization in super-resolution images as it describes the (regularized) optimal matching between the two protein intensity distributions. 
The set of pixels define the ground space $\mathcal{X}=\{x_1,\ldots,x_N\}$ with $N = N_x \cdot N_y$, where $N_x, N_y$ are the number of pixels in $x$- and $y$- direction, respectively. The pixels colour intensities are understood (after normalization) as discrete probability distributions $r,s \in \Delta_N$ on an equidistant grid in $[0,N_x\cdot l] \times [0,N_y \cdot l]$, where $l$ is the pixel size. The cost to transport pixel intensities from one pixel to the other is given by the squared euclidean distance $c_{(i-1)N+j}= \left\| x_i - x_j\right\|^2$ and $c_{\text{max}}\coloneqq \text{diam}(\mathcal{X})$ is the maximal distance on the ground space $\mathcal{X}$.
We introduce the \textit{regularized colocalization measure} {\small \textit{RCol}} which is based on the ROT plan $\pi_{p,\lambda,f}(r,s)$ and defined for $t \in [0,c_{\text{max}}]$ by 
\begin{equation}
\text{{\small RCol}}\coloneqq \text{{\small RCol}}(\pi_{p,\lambda,f}(r,s))(t) = \sum_{i=1}^{N^2} {\pi_{p,\lambda,f}(r,s)}_i \mathbbm{1}\{c_i\leq t\}\, .
\label{eq:rcol}
\end{equation}
Intuitively, $\text{{\small RCol}}(\pi_{p,\lambda,f}(r,s))(t)$ is the proportion of pixel intensities which is transported on scales smaller or equal to $t$ in the (regularized) optimal matching of the two intensity distributions with respect to some amount of regularization specified by $\lambda$. The function $\text{{\small RCol}}(\pi)$ constitutes an element in $\mathcal{D}[0,c_{\text{max}}]$ the space of all \textit{c\`{a}dl\`{a}g} functions \citep{billingsley2013convergence} on $[0,c_{\text{max}}]$ equipped with the supremum norm $\lVert f \rVert_{\infty} \coloneqq \sup\limits_{t \in [0,c_{\text{max}}]} \lvert f(t) \rvert$. 

\begin{theorem}\label{thm:colocalizationlimit}
Let $ \widehat{\text{{\small RCol}}}_n\coloneqq \text{{\small RCol}}(\pi_{p,\lambda,f}(\hat{r}_n,s))$ be the empirical regularized colocalization. Under the assumptions of \Cref{thm:optimalsolutionlimit}, as $n \to \infty$ it holds that
\begin{equation*}
\sqrt{n}\left\lbrace \widehat{\text{{\small RCol}}}_n-\text{{\small RCol}}\right\rbrace \overset{D}{\longrightarrow} \text{{\small RCol}}(G) \, ,
\end{equation*}
where $G$ is the centred Gaussian random variable in $\mathbb{R}^{N^2}$ with covariance $\Sigma_{p,\lambda,f}(r\vert s)$ given in \eqref{eq:covarianceregtransport}.
\end{theorem}

This provides \textit{approximate uniform confidence bands (CBs)} for the regularized colocalization. For $\alpha \in (0,1)$, let $\mathfrak{u}_{1-\alpha}$ be the $1-\alpha$ quantile from the distribution of $\lVert \text{{\small RCol}}(G)\rVert_{\infty}$. Hence,
\begin{equation}\label{eq:confidencebandone}
\mathcal{I}_n \coloneqq \left[ -\frac{\mathfrak{u}_{1-\alpha}}{\sqrt{n}} + \widehat{\text{{\small RCol}}}_n,\,\frac{\mathfrak{u}_{1-\alpha}}{\sqrt{n}} + \widehat{\text{{\small RCol}}}_n \right]
\end{equation}
is a $1-\alpha$ approximate uniform CB for the regularized colocalization $\text{{\small RCol}}$. More precisely, it holds that $\lim_{n \to \infty} \mathbb{P}(\text{{\small RCol}} \in \mathcal{I}_n) = 1-\alpha$.
In our subsequent data example we require to estimate both probability distributions $r,s\in \Delta_N$ (see \Cref{Estimationboth}). The confidence band \eqref{eq:confidencebandone} naturally extends to this case. Defining for $n,m\in \N$ the two sample empirical version of the regularized colocalization
\begin{equation*}
\widehat{\text{{\small RCol}}}_{n,m} \coloneqq \text{{\small RCol}}(\pi_{p,\lambda,f}(\hat{r}_n,\hat{s}_m))\, ,
\end{equation*}
we have, as the sample size $n \wedge m \to +\infty$ and $\nicefrac{m}{n+m}\to \delta \in (0,1)$, that 
\begin{equation}\label{eq:rcollimit}
\sqrt{\frac{nm}{n+m}}\left\lbrace \widehat{\text{{\small RCol}}}_{n,m}-\text{{\small RCol}}\right\rbrace \overset{D}{\longrightarrow} \text{{\small RCol}}(G) \, .
\end{equation}
Now, the random variable $G$ is centred Gaussian with covariance matrix $\Sigma_{p,\lambda,f}(r,s)$ given in \eqref{eq:twosamplecovariance}. In particular, for equal sample size $n=m$ the corresponding two sample CB reads as 
\begin{equation}\label{eq:confidencebandtwo}
\mathcal{I}_{n,n} \coloneqq \left[ -\frac{\sqrt{2}\mathfrak{u}_{1-\alpha}}{\sqrt{n}} + \widehat{\text{{\small RCol}}}_{n,n},\,\frac{\sqrt{2}\mathfrak{u}_{1-\alpha}}{\sqrt{n}} + \widehat{\text{{\small RCol}}}_{n,n} \right]\, ,
\end{equation}
with $\mathfrak{u}_{1-\alpha}$ according to the supremum of the r.h.s. of \eqref{eq:rcollimit}.

\subsection{Bootstrap Confidence Bands for Colocalization}

We denote by $\widehat{\text{{\small RCol}}}^*_{n,n}$ a bootstrap version of the regularized colocalization based on the empirical bootstrap estimators $\hat{r}_n^*$ and $\hat{s}_n^*$ derived by bootstrapping from the empirical distributions $\hat{r}_n$ and $\hat{s}_n$ (see \Cref{Bootstrap}). Note that $\text{{\small RCol}}$ is Lipschitz and therefore, according to \Cref{thm:bootstrap} and an application of the continuous mapping theorem for the bootstrap \citep[Proposition 10.7]{kosorok2008introduction}, we deduce that 
\begin{equation*}
\begin{split}
 \sup_{h \in \text{BL}_1(\mathbb{R})} \biggr\lvert &\mathbb{E}\left[ h\left(\left\lVert \sqrt{\nicefrac{n}{2}}\left\lbrace\widehat{\text{{\small RCol}}}^*_{n,n} - \widehat{\text{{\small RCol}}}_{n,n}\right\rbrace \right\rVert_{\infty}\right) \bigr\vert X_1,\ldots,X_n,Y_1,\ldots,Y_n \right] \\
 &-\mathbb{E}\left[ h\left(\left\lVert \sqrt{\nicefrac{n}{2}}\left\lbrace\widehat{\text{{\small RCol}}}_{n,n} - \text{{\small RCol}}\right\rbrace \right\rVert_{\infty}\right)\right]\biggr\rvert \overset{\mathbb{P}}{\longrightarrow} 0\, .
\end{split}
\end{equation*}
Hence, the quantile $\mathfrak{u}_{1-\alpha}$ is consistently approximated by its bootstrap analogue, say $\mathfrak{u}^*_{1-\alpha}$. This yields a bootstrap approximation of the CB in \eqref{eq:confidencebandtwo}. 

\subsubsection{Validation of Bootstrap and Resampling on Real Data}

\begin{figure}
\centering
\makebox{\includegraphics[width=0.68\textwidth]{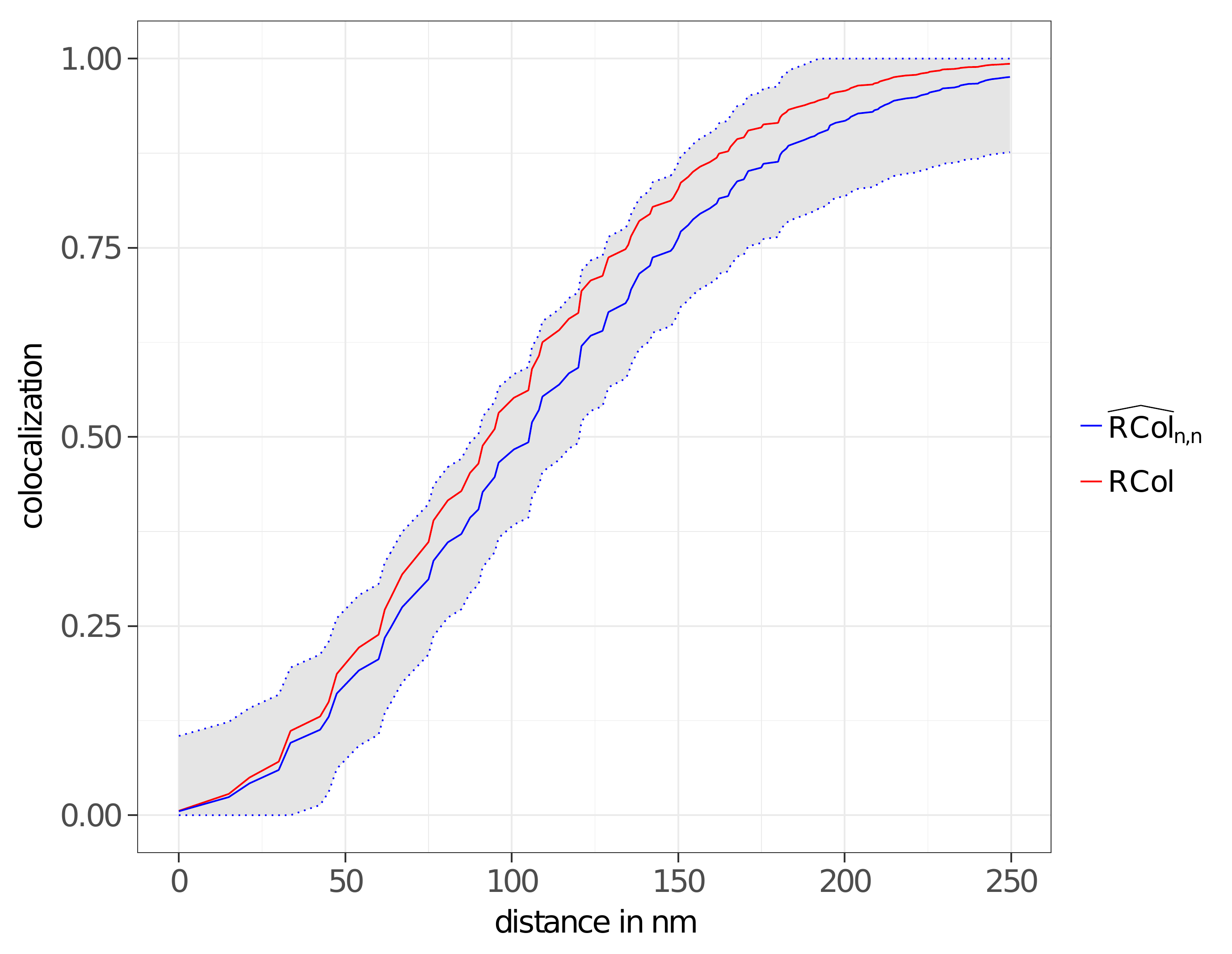}}
\caption{\label{fig:Regcolboot}\textbf{Entropy regularized colocalization for setting two: Staining of ATP Synthase and MIC60 for the zoom in ($\mathbf{128 \times 128}$ image).} The sampled regularized colocalization $(\lambda_0=0.01)$ (solid blue line, subsampling $n=2000$) with bootstrap confidence bands (gray area between dashed blue lines) based on the $n$ out of $n$ bootstrap with $B=100$ replications and $\alpha=0.05$. Red solid line: Population regularized colocalization.}
\end{figure}

The goal of the following analysis is to investigate the validity of the bootstrap CBs in combination with the resampling scheme (see \Cref{Subsampling}). To this end, we consider different pairs of zoom ins ($128\times 128$ sections) of the STED images in the first column in \Cref{fig:STED} including the pairs depicted in the middle column. For these instances we are still able to calculate the full regularized colocalization $\text{{\small RCol}}$ without resampling. Our goal is to validate that $\text{{\small RCol}}$ is covered by the bootstrap confidence bands. Statistically speaking after fixing a significance level $\alpha$, we are interested how close the empirical coverage probability is to the claimed nominal coverage of $1-\alpha$. The regularizer $f$ is the entropy \eqref{eq:entropy} with amount of regularization given by \eqref{eq:lambdaregularization} for specified $\lambda_0>0$. We resample $n=2000$ times according to the intensity distribution of the $128 \times 128$ image, calculate $B=100$ bootstrap replications and set the significance level for the CBs in \eqref{eq:confidencebandtwo} to $\alpha=0.05$. As an illustrative example, \Cref{fig:Regcolboot} demonstrates a case where $\text{{\small RCol}}$ is covered by the bootstrap CB. 

To investigate how well the empirical coverage probability approximates the nominal coverage probability of $1-\alpha=0.95$, we repeat our approach $100$ times and report how often $\text{{\small RCol}}$ is covered by the bootstrap CB. The result is given in \Cref{tab:Powerboot}. For bootstrap replications $B=100$ and rather large amount of regularization $\lambda_0\in\{0.5,\,1,\,2\}$ we are close to the nominal coverage probability. For small regularization $\lambda_0=0.01$ we obtain a slightly smaller empirical coverage probability than desired. This observation is consistent with our empirical simulations for the ROT distance in \Cref{Simulations}.

\begin{table}
\caption{\label{tab:Powerboot}\textbf{Validation of the bootstrap confidence bands.} Empirical coverage probability. Resampling $n=2000$, bootstrap replications $B=100$, regularization $\lambda_0$ and nominal coverage $1-\alpha = 0.95$.}
\centering
\fbox{\begin{tabular}{cc}
  Regularization $\lambda_0$  & Empirical coverage probability \\ \hline
   2 & 0.98 \\ 
   1 & 0.97 \\ 
   0.5 &  0.93 \\ 
   0.01 & 0.88 
   \end{tabular}}
\end{table}

\subsubsection{Empirical Colocalization Analysis of the STED Data} \label{subsec:HighResolution}

\begin{figure}
	\centering
	\includegraphics[width=0.68\textwidth]{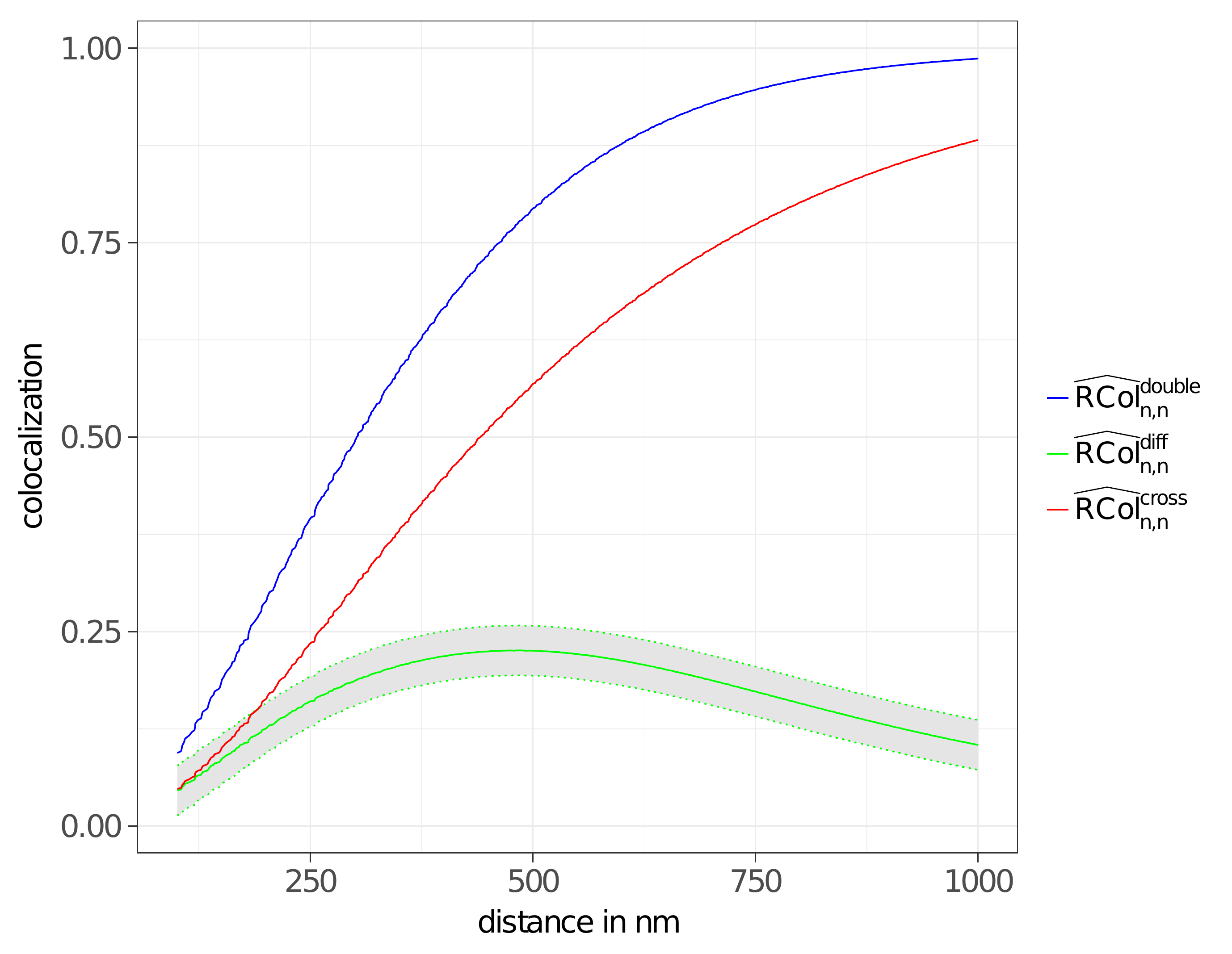}
	\caption{\textbf{Entropy regularized colocalization analysis.} Sampled regularized $(\lambda_0=0.01)$ colocalization (double staining of ATP Synthase, blue solid line; staining of ATP Synthase and MIC60, red solid line; subsampling $n=50.000$) based on the images on the left in Figure \ref{fig:STED} together with their difference (solid green line) with bootstrap confidence band (gray area between dashed green lines, $n$ out of $n$ bootstrap, $\alpha=0.05$, $B=100$ bootstrap replications).}	
	\label{fig:DiffBoot}
\end{figure}

We apply our resampling scheme on the full sized images ($666\times 666$ pixels) to evaluate the spatial proximity for each of the two settings, i.e. double staining of \textit{ATPS} (\Cref{fig:Regcolboot}) and staining of \textit{ATPS} and \textit{MIC60} (\Cref{fig:DiffBoot}). We expect the regularized colocalization for setting one (double staining of \textit{ATP Synthase}), say $\text{{\small RCol}}^{\text{double}}$, to be large in small distance regimes as most of the transport of pixel intensities should be carried out on small distances. For the regularized colocalization in setting two (staining of \textit{ATP Synthase} and \textit{MIC60}), say $\text{{\small RCol}}^{\text{cross}}$, we should observe that there is a significant amount of pixel intensities transported over larger distances resulting in a colocalization that is rather small in the small distance regimes.\\
Recall that in our validation setup we resampled $n=2000$ times out of pixel intensities represented by $128 \times 128$ images. To achieve comparable accuracy, we require here for $666\times 666$ pixel images a resampling scheme based on $n=50,000$ resamples of the underlying pixel intensity distributions. We then calculate their sampled regularized colocalization, denoted as $\widehat{\text{{\small RCol}}}_{n,n}^{\text{double}}$ and $\widehat{\text{{\small RCol}}}_{n,n}^{\text{cross}}$, respectively. In order to compare them, we propose to check their difference
$\widehat{\text{{\small RCol}}}_{n,n}^{\text{diff}} \coloneqq \widehat{\text{{\small RCol}}}_{n,n}^{\text{double}} - \widehat{\text{{\small RCol}}}_{n,n}^{\text{cross}}$,
especially in the small distance regime. As before, we obtain CBs by bootstrapping.\\
The results are presented in \Cref{fig:DiffBoot}. We observe that the sampled regularized colocalization $\widehat{\text{{\small RCol}}}_{n,n}^{\text{double}}$ (solid blue line) is larger than $\widehat{\text{{\small RCol}}}_{n,n}^{\text{cross}}$ (solid red line). Considering their difference $\widehat{\text{{\small RCol}}}_{n,n}^{\text{diff}}$ (solid green line) together with the bootstrap CB for the difference (gray area between dashed green lines, $\alpha=0.05$) demonstrates that the difference is significantly positive at all spatial scales below $1000$nm. In fact, our resampling approach for the regularized colocalization analysis reveals that double staining of \textit{ATP Synthase} is significantly more colocalized than staining of \textit{ATP Synthase} and \textit{MIC60} on all relevant spatial scales as biologically expected.

\section*{Acknowledgments}
The authors are grateful to S. Jakobs and S. Stoldt for providing us with the data in Section \ref{Colocalization} as well as to an AE and a referee for helpful comments.\\
Support by the Volkswagen Foundation and the DFG Research Training Group 2088 \textit{Discovering structure in complex data: Statistics meets Optimization and Inverse Problems} is gratefully acknowledged.

\bibliographystyle{plainnat.bst}
\bibliography{Sections/ROTlimit}{}

\appendix
\section{Proofs}\label{Appendix}

\subsection{Proof of Theorem \ref{thm:sensitivity}}\label{proof:sensitivity}
\textit{The first statement can be found in \citet{dessein2018regularized}. For the second statement we notice that for proper regularizer $f$ the ROT \eqref{eq:EntropyROT} with marginals $r_0$ and $s_0$ fulfils Slater's constraint qualification \citep[Proposition 26.18]{bauschke2017convex}. Hence, strong duality holds and the dual problem admits an optimal solution. Moreover, the ROT plan $\pi_{p,\lambda,f}$ and its corresponding optimal dual solution $\mu_{p,\lambda,f} \in \mathbb{R}^{2N-1}$ are characterized by the necessary and sufficient Karush-Kuhn-Tucker conditions
\begin{equation*}
c_p + \lambda [\nabla f(\pi_{p,\lambda,f})]^T - A_{\star}^T \mu_{p,\lambda,f} =0\, ,  \quad
A_{\star}\pi_{p,\lambda,f}-[r_0,{s_0}_{\star}]^T=0\, . 
\end{equation*}
The statement of the theorem now follows by an application of the implicit function theorem to this system of equalities. We define 
\begin{equation*}
G\colon \mathbb{R}^{N^2}\times \mathbb{R}^{2N-1} \times \mathbb{R}^{2N-1} \longrightarrow \mathbb{R}^{N^2+2N-1}
\end{equation*}
to be the function given by
\begin{equation*}
G(\pi,\mu,(r,s_{\star}))=\begin{bmatrix} c_p + \lambda [\nabla f(\pi)]^T - A_{\star}^T \mu \\ A_{\star}\pi-[r,\,s_{\star}]^T\end{bmatrix}\, ,
\end{equation*}
which is continuously differentiable in a neighbourhood of the specific point $(\pi_{p,\lambda,f}, \mu_{p,\lambda,f},(r_0,{s_0}_{\star}))$ with $G(\pi_{p,\lambda,f},\mu_{p,\lambda,f},(r_0,{s_0}_{\star}))=0$. The matrix of partial derivatives of $G$ with respect to $\pi$ and $\mu$ is given by
\begin{align*}
\nabla_{(\pi,\mu)}G(\pi_{p,\lambda,f},\mu_{\lambda,f},&(r_0,{s_0}_{\star})) 
= \\
&\begin{pmatrix}
\lambda \nabla^2 f(\pi_{p,\lambda,f}) & -A_{\star}^T \\
A_{\star} & 0
\end{pmatrix} \in \mathbb{R}^{N^2+2N-1\times N^2+2N-1}\, .
\end{align*}
and is non-singular as $\lambda >0$, the Hessian $\nabla^2 f(\pi_{p,\lambda,f})$ is positive definite by definition of a proper regularizer and the matrix $A_{\star}^T$ has full rank. The implicit function theorem guarantees the existence of a neighbourhood $\mathcal{U}$ around $(r_0,{s_0}_{\star})$ and a continuously differentiable function $\Theta\colon \mathcal{U}\to \mathbb{R}^{N^2}\times \mathbb{R}^{2N-1}$ with components $\Theta(r,s_{\star})=(\phi_{p,\lambda,f}(r,s_{\star}),\gamma_{p,\lambda,f}(r,s_{\star}))$ such that $$G(\phi_{p,\lambda,f}(r,s_{\star}),\gamma_{p,\lambda,f}(r,s_{\star}),(r,s_{\star}))=0$$ for all $(r,s_{\star})\in \mathcal{U}$. The vector parametrized by $\phi_{p,\lambda,f}(r,s_{\star})$ fulfils the necessary and sufficient optimality conditions for the ROT with $(r,s_{\star})\in \mathcal{U}$. Hence, the function $\phi_{p,\lambda,f}\colon \mathcal{U}\to \mathbb{R}^{N^2}$ parametrizes the ROT plan with feasible set $\mathcal{F}(r,s)$ for all $(r,s_{\star})\in \mathcal{U}$.
It remains to compute the derivative of $\phi_{p,\lambda,f}$ at $(r_0,{s_0}_{\star})$. According to the implicit function theorem we obtain that 
\begin{equation*}
\begin{split}
\nabla \Theta (r_0,{s_0}_{\star}) &= \begin{bmatrix}
\nabla \phi_{p,\lambda,f}(r_0,{s_0}_{\star})\\ \nabla \gamma_{p,\lambda,f}(r_0,{s_0}_{\star})
\end{bmatrix}\\
&=[\nabla_{(\pi,u)} G(\pi_{p,\lambda,f},\mu_{p,\lambda,f},(r_0,{s_0}_{\star}))]^{-1} \begin{bmatrix}
0_{N^2\times 2N-1} \\ I_{2N-1\times 2N-1}
\end{bmatrix}\\
&=\begin{bmatrix}
[\nabla^2 f(\pi_{p,\lambda,f})]^{-1}\,A_{\star}^T\,[A\,[\nabla^2 f(\pi_{p,\lambda,f})]^{-1}\,A_{\star}^T]^{-1}\\
-[A_{\star}\,[\nabla^2f(\pi_{p,\lambda,f})]^{-1}\,A_{\star}^T]^{-1}
\end{bmatrix}\, ,
\end{split}
\end{equation*}
where in the second equality the indices denote the size of the zero and identity matrix, respectively, For the last equality we applied a result by \citet[equalities 4.2.8]{fiacco1984sensitivity}. \qed }

\subsection{Proof of Theorem \ref{thm:optimalsolutionlimit}}\label{proof:optimalsolutionlimit}
\textit{Since the vector $n\hat{r}_n$ follows a multinomial distribution with probability vector $r$, the multivariate central limit theorem yields
\begin{align*}
\sqrt{n}(\hat{r}_n - r) \overset{D}{\longrightarrow} \mathcal{N}_N (0,\Sigma(r))\, .
\end{align*}
The multivariate delta method in conjunction with the sensitivity analysis for regularized transport plans concludes the statement. More precisely, we obtain that
\begin{align*}
\sqrt{n}\left\lbrace \pi_{\lambda,f}(\hat{r}_n,s)-\pi_{\lambda,f}(r,s) \right\rbrace 
= \sqrt{n}\left\lbrace \phi_{p,\lambda,f}(\hat{r}_n,s_{\star})-\phi_{p,\lambda,f}(r,s_{\star}) \right\rbrace\\
 \overset{D}{\longrightarrow} \mathcal{N}_{N^2}(0,\nabla_r\, \phi_{p,\lambda,f}(r,s_{\star}) \, \Sigma(r)\, [\nabla_r\, \phi_{p,\lambda,f}(r,s_{\star})]^T) \, .
\end{align*}
This finishes the proof. \qed}

\subsection{Proof of Theorem \ref{thm:optimaltransportlimit}}\label{proof:optimaltransportlimit}
\textit{Let $G \sim \mathcal{N}_{N^2}(0,\Sigma_{\lambda,f}(r \vert s)) $ be a $N^2$-dimensional Gaussian random vector distributed according to the limit distribution in Theorem \ref{thm:optimalsolutionlimit}. The multivariate delta method yields
\begin{align*}
&\sqrt{n}\left\lbrace W_{p,\lambda,f} (\hat{r}_n,s)-W_{p,\lambda,f}(r,s) \right\rbrace  \\
=&\sqrt{n}\, \left\lbrace\langle c, \pi_{p,\lambda, f}(\hat{r}_n,s) \rangle^{\frac{1}{p}}- \langle c, \pi_{p,\lambda, f}(r,s) \rangle^{\frac{1}{p}} \right\rbrace
 \overset{D}{\longrightarrow} \gamma^T \,G\, ,
\end{align*}
where $\gamma$ is the gradient of the function $\pi \mapsto \langle c, \pi \rangle^{\frac{1}{p}}$ evaluated at the regularized transport plan $\pi_{p,\lambda,f}(r,s)$. The variance of the real-valued random variable $\gamma^T \,G$ is given by
\begin{align*}
\sigma^2_{p,\lambda,f}(r\vert s)=\gamma^T\, \Sigma_{p,\lambda,f}(r \vert s) \, \gamma \, .
\end{align*}
In particular, the variance $\Sigma_{p,\lambda,f}(r \vert s)$ is continuous in $r$. Consequently, by the strong law of large numbers $\hat{r}_n\overset{a.s.}{\longrightarrow} r$ and the continuous mapping theorem we have that $\sigma^2_{p,\lambda,f}(\hat{r}_n \vert s) \overset{a.s.}{\longrightarrow} \sigma^2_{p,\lambda,f}(r\vert s)$ which together with Slutzky's theorem \citep[Lemma 2.8]{van2000asymptotic} concludes the statement. \qed}

\subsection{Proof of Theorem \ref{thm:colocalizationlimit}}\label{proof:colocalizationlimit}
\textit{The map $\pi \mapsto \text{RCol}(\pi)$ is linear and $1$-Lipschitz. Hence, according to our distributional limit results in Theorem \ref{thm:optimalsolutionlimit} for the regularized optimal transport plan and the continuous mapping \citep[Theorem 18.11]{van2000asymptotic} we conclude that 
\begin{align*}
\sqrt{n}\left\lbrace \widehat{\text{RCol}}_n-\text{RCol}\right\rbrace 
= \text{RCol}(\sqrt{n}\left\lbrace \pi_{p,\lambda,f}(\hat{r}_n,s)-\pi_{p,\lambda,f}(r,s)\right\rbrace) \overset{D}{\longrightarrow} \text{RCol}(G)\, ,
\end{align*}
where $G\sim \mathcal{N}_{N^2}(0,\Sigma_{\lambda,f}(r \vert s))$. \qed}

\end{document}